\documentclass[11pt]{article}

\usepackage{fullwidth}
\usepackage{color}
\usepackage[color,matrix,arrow]{xy}
\usepackage[top=1.5in, bottom=1.5in, left=1in, right=1in]{geometry}
\usepackage{amsgen}
\usepackage{amsmath}
\usepackage{amstext}
\usepackage{amsbsy}
\usepackage{amsopn}
\usepackage{amsfonts}
\usepackage{amssymb}
\usepackage{eepic}
\usepackage{pgf,tikz}
\usepackage{graphicx}
\usepackage{epsf}
\usepackage{pstricks}
\xyoption{all}

\def\Box{\square}
\def\edge{\relbar\joinrel\relbar}

\def\tra#1{\smash{\mathop{\mid\kern
-1pt\joinrel\relbar\joinrel\relbar}\limits^{*}_{#1}}}
\def\longtra#1{\smash{\mathop{\mid\kern
-1pt\joinrel\relbar\joinrel\relbar\joinrel\relbar}\limits^{*}_{#1}}}
\def\vlongtra#1{\smash{\mathop{\mid\kern
-1pt\joinrel\relbar\joinrel\relbar\joinrel\relbar\joinrel\relbar}\limits^{*}_{#1}}}
\def\vvlongtra#1{\smash{\mathop{\mid\kern
-1pt\joinrel\relbar\joinrel\relbar\joinrel\relbar\joinrel\relbar\joinrel\relbar}\limits^{*}_{#1}}}
\def\vvvlongtra#1{\smash{\mathop{\mid\kern
-1pt\joinrel\relbar\joinrel\relbar\joinrel\relbar\joinrel\relbar\joinrel\relbar\joinrel\relbar}\limits^{*}_{#1}}}
\def\etra#1{\smash{\mathop{\mid\kern
-1pt\joinrel\relbar\joinrel\relbar}\limits_{#1}}}

\def\on#1{\{1,\ldots,{#1}\}}

\def\Rw{\Rightarrow}
\def\oo{\overline}

\def\wh{\widehat}
\def\B{{\cal{B}}}

\def\F{{\cal{F}}}

\def\L{{\cal{L}}}

\def\S{{\cal{S}}}

\def\fct{\mbox{fct}}

\def\dim{\mbox{dim}}

\def\per{\mbox{Per}}

\def\tbpav{\mbox{TBPav}}
\def\bpav{\mbox{BPav}}
\def\pav{\mbox{Pav}}

\def\ker{\mbox{Ker}\,}

\def\max{\mbox{max}}

\def\cl{\mbox{Cl}}
\def\rk{\mbox{rk}}

\def\flatx{\mbox{Fl}}

\def\H{{\mathcal{H}}}

\def\V{{\cal{V}}}

\def\Z{\mathbb{Z}}

\def\p{\varphi}

\def\inv{^{-1}}

\def\tbr{{\cal{TBR}}}
\def\tbp{{\cal{TBP}}}
\def\pure{\mbox{pure}}
\def\SB{\mathbb{SB}}
\def\bi{\begin{itemize}}
\def\ei{\end{itemize}}
\def\beq{\begin{equation}}
\def\eeq{\end{equation}}

\def\J{{\cal{J}}}
\def\Y{{\cal{Y}}}

\newtheorem{T}{Theorem}[section]
\newcommand{\bt}{\begin{T}}
\newcommand{\et}{\end{T}}
\newcommand{\ftd}{$\square$\end{T}}

\newtheorem{Proposition}[T]{Proposition}
\newcommand{\bp}{\begin{Proposition}}
\newcommand{\ep}{\end{Proposition}}
\newcommand{\fpd}{$\square$\end{Proposition}}

\newtheorem{Lemma}[T]{Lemma}
\newcommand{\bl}{\begin{Lemma}}
\newcommand{\el}{\end{Lemma}}
\newcommand{\fld}{$\square$\end{Lemma}}

\newtheorem{Corol}[T]{Corollary}
\newcommand{\bc}{\begin{Corol}}
\newcommand{\ec}{\end{Corol}}
\newcommand{\fcd}{$\square$\end{Corol}}

\newtheorem{Result}[T]{Result}
\newcommand{\br}{\begin{Result}}
\newcommand{\er}{\end{Result}}
\newcommand{\frd}{$\square$\end{Result}}

\newtheorem{Example}[T]{Example}
\newcommand{\be}{\begin{Example}}
\newcommand{\ee}{\end{Example}}

\newtheorem{Problem}[T]{Problem}
\newcommand{\bq}{\begin{Problem}}
\newcommand{\eq}{\end{Problem}}

\newtheorem{Remark}[T]{Remark}
\newcommand{\brem}{\begin{Remark}}
\newcommand{\erem}{\end{Remark}}

\newtheorem{Conj}[T]{Conjecture}
\newcommand{\bj}{\begin{Conj}}
\newcommand{\ej}{\end{Conj}}
\newcommand{\fjd}{$\square$\end{Conj}}

\newcommand{\proof}
   {\par\medbreak\noindent{\bf Proof}.\enspace}

\newcommand{\qed}{%\hfill
$\Box$
\par\bigbreak}

\normalbaselines
\def\abstract#1{\par\bigskip
\begingroup\small
\baselineskip=12truept
\begin{center}ABSTRACT\end{center}
\par\medskip\par\noindent
\null\hfill\hbox{\vbox{\hsize=5truein\noindent#1}}
\hfill\null\par\endgroup\par}

\title{Truncated boolean representable simplicial complexes}
\author{{\bf Stuart Margolis, John Rhodes and Pedro V. Silva}}

\date{\today}

\begin{document}
\maketitle
% ----------------------------------------------------------------

\begin{center}\small

\bigskip

\end{center}

\abstract{We extend, in significant ways, the brief theory of truncated boolean representable simplicial complexes introduced in 2015. This theory, which includes all matroids, represents the largest class of finite simplicial complexes for which combinatorial geometry can be meaningfully applied.
}

\section{Introduction}

In this paper, we extend  the theory of TBRSC (truncated boolean representable simplicial complexes) created in \cite[Sec. 8.2]{RSm}. The paper is reasonably self-contained, but familiarity with \cite{RSm} will be very useful. 

Matroids \cite{Oxl}, BRSC (boolean representable simplicial complexes) and TBRSC, as models of {\em discrete geometry}, are all concerned with the generalized notion of {\em independence}. All matroids admit a boolean representation (usually many), so do BRSC, but not conversely, so BRSC are beyond matroids \cite{Cam}. The set of {\em independent} subsets of a finite set of points $V$ form a simplicial complex $(V,H)$ (in the sense of elementary algebraic topology \cite{Spa}), $H$ being a nonempty collection of subsets of $V$ closed under taking subsets.

We are interested in simplicial complexes $S = (V,\H)$ {\em arising from a geometry}, as to be explained below. In this introduction we only concern ourselves with simple simplicial complexes, i.e. those such that all pairs $v_1v_2$ of (distinct) elements of $V$ are in $\H$. 
But non-simple simplicial complexes are also considered in \cite{RSm} and several papers. 

Simple matroids arise through transversals of the partial differences for chains in geometric lattices, where we identify the vertices with the atoms of the lattice. More generally, a BRSC may be obtained through transversals of the partial differences for chains in an arbitrary lattice, the role of atoms being played by any join-generating set \cite[Chapter 5]{RSm}.

Given a simplicial complex $S = (V,\H)$, the combinatorial and algebraic fields use the term {\em rank} $r$ as the cardinality of the largest set in $\H$. The topological and geometric fields use $d = r-1$ as the {\em dimension} of $S$. We use both as will be explained. We say that $S$ is {\em paving} (of dimension $d$) and write $S \in {\rm Pav}(d)$ if $S$ has dimension $d$ and $\H$ contains all subsets of $V$ of cardinality $d$. Thus we may identify the class of (finite) graphs with ${\rm Pav}(1)$. 

For (T)BRSC the geometry {\em comes in} similarly to matroids as will be explained. Let us restrict to the paving case for simplicity. Let $M = (V,\H)$ be a matroid in Pav$(d)$. For simplicity of explanation we assume that $d = 2$, but all generalizes to arbitrary $d$.

The lattice of flats $L(M)$ of $M$ induces a closure operator on $2^V$ (see \cite[Section 4.2]{RSm} for the most general version). Let $\oo{X}$ denote the closure of $X \subseteq V$. Let 
$$\L = \{ F \in L(M) \mid 2 \leq |F| < |V| \}.$$
If $M$ is a matroid, then each pair of points of $V$ is contained in a unique block in $\L$. 
This makes $\L$ a PBD ({\em pairwise balanced design}) (with $\lambda = 1$) in {\em design theory}. Conversely, and except for trivial cases, every such PBD determines a matroid in Pav$(2)$ (see e.g \cite[Proposition 4.2]{MRS2}).

This generalizes to $\bpav(2)$ (the BRSC in $\pav(2)$) as follows (see \cite[Section 6.3]{RSm}). Let $S = (V,\H) \in \bpav(2)$. Then 
$\L = \{ F \in L(S) \mid 2 \leq |F| < |V| \}$  is a {\em partial geometry} since $|F_1 \cap F_2| \leq 1$ for all distinct $F_1,F_2 \in \L$. This implies that each pair of points of $V$ is contained in at most one block in $\L$ (pretty much the central notion). If $S$ is not a matroid, this is not enough to produce a PBD, but we get a PBD by adding to $\L$ all pairs of points of $V$ contained in no element of $\L$.

If $S = (V,\H)$ is a simplicial complex and $k \geq 1$, then $T_k(S) =  (V,\H \cap P_{\leq k}(V))$ is the {\em truncation} of $S$ to rank $k$. Note that we use rank here, not dimension.

In general, the truncation of a BRSC is not a BRSC (so we get a wider concept denoted by TBRSC), but it is easily characterized as follows.

Given a a simplicial complex $S = (V,\H)$ of dimension $d$, we define
$$\varepsilon(S) = \varepsilon(\H) =  
\{ X \subseteq V \mid \forall Y \in \H \cap P_{\leq d}(X)\; \forall p \in V
\setminus X \hspace{.3cm} Y \cup \{ p \} \in \H\}.$$
Then transversals of the partial differences for chains in $\varepsilon(S)$ define a BRSC denoted by $S^{\varepsilon}$.

In the paving case, we actually have $\varepsilon(S) = L(S^{\varepsilon})$. Moreover, $S^{\varepsilon}$ truncated to the rank of $S$ is the unique largest TBRSC $S^0$ contained in $S$. We mean largest with respect to inclusion of faces (for the same vertex set). This relation is called the {\em weak order} in matroid theory \cite[Section 7.3]{Oxl}.

Now $S$ is a TBRSC if and only if $T_{d+1}(S^{\varepsilon}) = S$. This provides a useful criterion to recognize a TBRSC. 

We note that even if $M$ is a matroid, $M^{\varepsilon}$ is not necessarily a matroid (see \cite[Examples 5.5 and 5.6]{MRS2}). Indeed, $\varepsilon(M)$ constitutes a lattice under intersection, but not necessarily a geometric lattice. See \cite{MRS2}, where this is further developed. 

We note that, if $M$ is a matroid of rank $r$, then $\varepsilon(M)$ consists precisely of the $(r-1)$-{\em closed} subsets of $V$ in the terminology introduced by Crapo \cite{Cra}.
But TBRSC also shed a new light on Crapo's concept of erection. An erection of a matroid $M = (V,\H)$ of rank $r$ is a matroid $M'$ of rank $\leq r+1$ such that $M = T_r(M')$. In \cite{Cra}, Crapo proved that a collection $\{ B_1,\ldots, B_k\} \subset 2^V$ form a collection of maximal flats of an erection of $M$ if and only if:
\begin{itemize}
\item[(1)]
the closure of each $B_i$ in $L(M)$ is $V$;
\item[(2)]
each $B_i$ is $(r-1)$-closed (i.e. belongs to $\varepsilon(M)$);
\item[(3)]
every facet of $M$ (i.e. maximal subset in $\H$) is contained in a unique $B_i$.
\end{itemize}
He then proves that the collection of all erections of $M$ (including the trivial erection $M$) form a lattice for the weak order, and the maximal erection  is called the free erection.

These conditions were designed to remain in the matroid context. We go beyond matroids to the wider class of TBRSC and consider arbitrary differences of rank, generalizing the work of Crapo and others \cite{Cra, Ngu,PP}. For a TBRSC $S$, these operations can be viewed as {\em strong maps} ($\vee$-maps) of $S^{\varepsilon}$ (see \cite{MRS} and Chapter 5, especially Sections 5.4 - 5.5 and 8.2 of \cite{RSm}).

The outline of this paper is as follows:

In Sections 2 and 3 we provide the basic theory of BRSC and TBRSC, respectively. In Section 4 we discuss low dimensions.

Section 5 deals with the join operator: $(V,\H) \cup (V,\H') = (V,\H\cup \H')$. 
In this paper we use join referring to the lattice of all simplicial complexes on a fixed vertex set $V$, ordered by inclusion. This is the same as the lattice of semigroup ideals of the monoid $(2^V, \cap)$.

In general, (T)BRSC are not closed under join, but the class $\tbpav(d)$ (the TBRSC in $\pav(d)$) is. A key resource are the complexes $B_d(V,L)$ (where $2 \leq d \leq |L| < |V|$), containing all subsets of $V$ with at most $d$ points and all subsets of $V$ with $d+1$ points which intersect $L$ in exactly $d$ points. 

Note that, given $S = (V,\H) \in \tbpav(d)$, we have $B_d(V,L) \subseteq S$ if and only if $L \in \varepsilon(S)$. Moreover, $S$ is a TBRSC if and only if 
$$S = \bigcup_{L \in \L} B_d(V,L) = T_{d+1}(S^{\varepsilon})$$
where $\L \subseteq P_{\geq d}(V) \setminus \{ V \}$. More generally, $S^0 = T_{d+1}(S^{\varepsilon})$ is the largest paving TBRSC contained in $S$, and is therefore the largest subcomplex of $S$ allowing some geometrical features. This strictly includes all paving matroids.

In Section 6, we show that the maximum number of vertices for a minimal $S \in {\rm TBPav}(d) \setminus {\rm BPav}(d)$ is $(d+1)(d+2)$. On the other hand, the prevariety ${\cal{TBP}}_2$ (consisting of all paving TBRSC of dimension $\leq 2$) is not finitely based.

In Section 7, we discuss three questions involving the largest pure subcomplex of a BRSC or of one of its truncations. We answer them negatively in the general case, but we show them to hold for low dimensions,

Finally, we discuss in Section 8 some of the topological properties of the geometric realization of a TBRSC, generalizing previous results for BRSC.

\section{Boolean representable simplicial complexes}

For the material presented in this section, the reader is referred to \cite{RSm}. All the results mentioned here will be used throughout the paper without further reference.

All lattices and simplicial complexes in this paper are assumed to be
finite. Given a set $V$ and $n \geq 0$, we denote by $P_n(V)$
(respectively $P_{\leq n}(V), P_{\geq n}(V)$) the set
of all subsets of $V$ with precisely (respectively at most, at least) $n$
elements. 

A (finite) simplicial complex is a structure of the form $S =
(V,\H)$, where $V$ is a finite nonempty set and $\H \subseteq 2^V$
contains $P_1(V)$ and is closed under taking subsets. 
The elements of $V$ and $\H$ are called respectively {\em vertices} and {\em faces}. 
To
simplify notation, we shall often denote a face $\{
x_1,x_2,\ldots,x_n\}$ by $x_1x_2\ldots x_n$. 

A face of $S$ which is maximal with respect to
inclusion is called a {\em facet}. We denote by ${\rm fct}(S)$ the set
of facets of $S$.
The rank and dimension of $S$ are defined respectively by
$$\rk(S) = \max\{ |I| : I \in \H\},\quad 
\dim(S) = \rk(S) -1.$$

We say that $S = (V,\H)$ is:
\bi
\item {\em simple} if $P_2(V) \subseteq \H$;
\item {\em paving} if $P_{\dim(S)}(V) \subseteq V$.
\ei
We denote by $\pav(d)$ the class of all paving simplicial complexes of dimension $d$.

Two simplicial complexes $(V,\H)$ and $(V',\H')$ are {\em isomorphic} if
 there exists a bijection $\p:V \to V'$ such that
 $$X \in \H \mbox{ if and only if }X\p \in \H'$$
 holds for every $X \subseteq V$.

 If $S = (V,\H)$ is a simplicial complex and 
 $W \subseteq V$ is
 nonempty, we call 
 $$S|_W = (W,\H \cap 2^W)$$
 the {\em restriction} of $S$ to $W$. It is obvious that $S|_W$ is
 still a simplicial complex. 

A simplicial complex $M = (V,\H)$ is called a {\em matroid} if it
satisfies the {\em exchange property}:
\bi
\item[(EP)]
For all $I,J \in \H$ with $|I| = |J|+1$, there exists some
  $i \in I\setminus J$ such that $J \cup \{ i \} \in \H$.
\ei

An important example of matroids are the {\em uniform matroids} $U_{k,n}$: for all $1 \leq k \leq n$, we write $U_{k,n} = (V,P_{\leq k}(V))$ with $|V| = n$.

Given an $R \times V$ matrix $M$ and $Y \subseteq R$, $X \subseteq V$,
we denote by $M[Y,X]$ the submatrix of $M$ obtained by deleting all
rows (respectively columns) of $M$ which are not in $Y$ (respectively $X$).

A boolean matrix $M$ is {\em lower unitriangular} if it is of the form
$$\left(
\begin{matrix}
1&&0&&0&&\ldots&&0\\
?&&1&&0&&\ldots&&0\\
?&&?&&1&&\ldots&&0\\
\vdots&&\vdots&&\vdots&&\ddots&&\vdots\\
?&&?&&?&&\ldots&&1
\end{matrix}
\right)
$$

Two matrices are {\em congruent} if we can transform one into the other
by independently permuting rows/columns. A boolean matrix is {\em nonsingular}
if it is congruent to a lower unitriangular matrix. 

Equivalently, nonsingular matrices can be characterized through the concept of {\em permanent}. 
The permanent of a square matrix $M = (m_{ij})$ (a 
positive version of the determinant) is defined by 
$$\per M = \displaystyle\sum_{\pi \in S_n} \prod_{i=1}^n
m_{i,i\pi}.$$
But, even though our matrix is boolean, we compute its permanent in the {\em superboolean} semiring $\SB$, which can be described as the quotient of the usual semiring $(\mathbb{N},+,\cdot)$ by the congruence with classes $\{ 0 \}, \{ 1 \}, \{ 2,3 \ldots \}$. Then a square boolean matrix is nonsingular if and only if its permanent is 1 in $\SB$.

We note that the
classical results on determinants involving only a rearrangement of
the permutations extend naturally to permanents.

Given an $R \times V$ boolean matrix $M$, we say that the
subset of columns $X \subseteq V$ is $M$-{\em independent} if there exists
some $Y \subseteq R$ such that the submatrix $M[Y,X]$ is nonsingular.

A simplicial complex $S = (V,\H)$ is {\em boolean representable} (BRSC) if
there exists some boolean matrix $M$ such that
$\H$ is the set of all $M$-independent subsets of $V$. Since $P_1(V) \subseteq \H$, this implies that all the columns of $M$ are nonzero. Moreover, for all distinct $p,q \in V$, the columns $M[R,p]$ and $M[R,q]$ are different if and only if $pq \in \H$.

By restricting the set of columns, it is easy to see that a restriction of a BRSC is still a BRSC \cite[Proposition 8.3.1(i)]{RSm}.

All matroids are boolean representable \cite[Theorem
  5.2.10]{RSm}, but the converse is not true.
  
A subset $\F$ of $2^V$ is called a {\em Moore family} if $V \in \F$ and $\F$ is closed under intersection (that is, a Moore family is a submonoid of the monoid of all subsets of $V$ under intersection). Every Moore family, under inclusion, constitutes a lattice (with intersection as meet and the determined join $F_1 \vee F_2 = \cap\{ F \in \F \mid F_1 \cup F_2 \subseteq F\}$). We say that $X \subseteq V$ is a {\em transversal of the
successive differences} for a chain
$$F_0 \subset F_1 \subset \ldots \subset F_k$$
in $\F$ if $X$ admits an enumeration $x_1,\ldots , x_k$ such that $x_i \in F_i
\setminus F_{i-1}$ for $i = 1,\ldots,k$. 

If ${\rm Tr}(\F)$ is the set of transversals of the
successive differences for chains in $\F$, then $(V,{\rm Tr}(\F))$ constitutes a BRSC. Moreover, every BRSC can be obtained this way by taking as Moore family its {\em lattice of flats} (see \cite[Chapters 5 and 6]{RSm}):

We say that $X
\subseteq V$ is a {\em flat} of $S = (V,\H)$ if
$$\forall I \in \H \cap 2^X \hspace{.2cm} \forall p \in V \setminus X
\hspace{.5cm} I \cup \{ p \} \in \H.$$
The set of all flats of $S$ is denoted by 
$L(S)$. Note that $V, \emptyset \in L(S)$ in all cases, and $L(S)$ is indeed a Moore family (and therefore a lattice).
Note also that $P_{\leq d-1}(V) \subseteq L(S)$ for every $S \in \pav(d)$. 

It follows from \cite[Corollary 5.2.7]{RSm} that a simplicial complex $S = (V,\H)$ is boolean
representable If and only if the transversals of the successive differences for chains in $L(S)$ are precisely the elements of $\H$. 

By \cite[Proposition 8.3.3(i)]{RSm}, the flats of a BRSC determine flats on any restriction: if $F$ is a flat of a BRSC $S = (V,\H)$ and $W \subseteq V$, then $F \cap W \in L(S|_W)$. 

The lattice $L(S)$ induces a closure operator on $2^V$ defined by
$$\oo{X} = \cap\{ F \in L(S) \mid X \subseteq F \}$$
for every $X \subseteq V$. 
It follows from the definitions that $\oo{X} = V$ when $X$ contains a facet of $S$.

By \cite[Corollary 5.2.7]{RSm}, $S = (V,\H)$ is boolean
representable if and only if every $X \in \H$ admits an enumeration
$x_1,\ldots, x_k$ satisfying
\beq
\label{derby3}
\oo{x_1} \subset \oo{x_1x_2} \subset\ldots \subset \oo{x_1\ldots x_k}.
\eeq
It is well known that in  the case of matroids, this enumeration can be chosen arbitrarily \cite{Oxl}.

\section{Truncation}

In this section, we exposit the basic facts about TBRSCs. The proofs of the results can be found in \cite[Section 8.2]{RSm}, but we include them in the Appendix for the sake of completeness and convenience for the reader.

Given a simplicial complex $S = (V,\H)$ and $k \geq 1$, the $k$-{\em
  truncation} of $S$ is the
simplicial complex $T_k(S) = (V,T_k(\H))$, where $T_k(\H) = \H \cap P_{\leq k}(V)$. 

We say that a simplicial complex $S = (V,\H)$ is a {\em TBRSC} if $S = T_k(S')$ for some BRSC $S'$ and $k \geq 1$. For every $d \geq 1$, we denote by $\tbpav(d)$ the class of all paving TBRSCs of dimension $d$.

To recognize a TBRSC, it is convenient to develop an alternative characterization. The key is building the flats of a canonical BRSC having our TBRSC as a truncation. The following result characterizes the flats of a truncation with respect to the flats of the original complex.

\bp
\label{tru}
{\rm \cite[Proposition 8.2.2]{RSm}}
Let  $S = (V,\H)$ be a simplicial complex and let $k \geq 1$. Then
$$L(T_k(S)) = \{ X \in L(S) \mid {\rm fct}(T_k(S)) \cap\,
2^X = \emptyset\} \cup \{ V \}.$$
\ep

\proof
In the Appendix.
\qed

It follows that the lattice of flats of $T_k(S)$ is obtained from the lattice of flats of $S$ by identifying the elements of an up set (namely the subset of flats containing some facet of $T_k(S)$). In semigroup-theoretic terms, this makes $L(T_k(S))$ a Rees quotient of the $\vee$-semilattice of $L(S)$.

For any simplicial complex $S = (V,\H)$ of dimension $d$, we define
$$\varepsilon(S) = \varepsilon(\H) = 
\{ X \subseteq V \mid \forall Y \in \H \cap P_{\leq d}(X)\; \forall p \in V
\setminus X \hspace{.3cm} Y \cup \{ p \} \in \H\}.$$

Note that $\varepsilon(S)$ generalizes to arbitrary simplicial complexes what Crapo calls {\em d-closed sets} of a monoid in his fundamental paper from 1970 \cite{Cra}.

The following lemma is clear from the definition.

\bl
\label{propt}
{\rm \cite[Lemma 8.2.3]{RSm}}
Let $S = (V,\H)$ be a simplicial complex. Then:
\bi
\item[(i)] $\varepsilon(S)$ is a Moore family;
\item[(ii)] $L(S) \subseteq \varepsilon(S)$.
\ei
\el

Given $S = (V,\H)$, write $\H^{\varepsilon} = {\rm Tr}(\varepsilon(S))$ and 
let $S^{\varepsilon} = (V,\H^{\varepsilon})$ denote the BRSC defined by the lattice $\varepsilon(S)$.

\bl
\label{inc}
{\rm \cite[Lemma 8.2.4]{RSm}}
Let $S = (V,\H)$ be a simplicial complex of dimension $d$. Then: 
\bi
\item[(i)] $T_{d+1}(\H^{\varepsilon}) \subseteq \H$;
\item[(ii)] $\varepsilon(S) \subseteq L(S^{\varepsilon})$;
\item[(iii)] $S^{\varepsilon}$ is a BRSC.
\ei
\el

\proof
In the Appendix.
\qed

Now we can state the main result of this section:

\bt
\label{eqtr}
{\rm \cite[Theorem 8.2.5]{RSm}}
Let $S = (V,\H)$ be a simplicial complex of dimension $d$. Then the following
conditions are equivalent:
\bi
\item[(i)] $S = T_{d+1}(S')$ for some boolean representable simplicial
  complex $S'$; 
\item[(ii)] $S = T_{d+1}(S^{\varepsilon})$.
\ei
Furthermore, in this case we have $L(S^{\varepsilon}) = \varepsilon(S)$.
\et

\proof
In the Appendix.
\qed

We present now two examples which show that the class of TBRSCs is intermediate between the class of BRSCs and the class of simplicial complexes. We analyze these examples in the Appendix.

The first example shows that a TBRSC is not necessarily a BRSC, even in the paving case.

\be
\label{truno}
Let $V = \on{6}$, 
$\H = (P_{\leq 3}(V) \setminus \{ 135, 235,
146, 246, 346, 456 \})$
and $S = (V,\H)$. Then $S \in {\rm TBPav}(2)$
but is not boolean representable.
\ee

The second complex shows that a (paving) simplicial complex is not necessarily a TBRSC.

\be
\label{extruA}
Let $V = \on{4}$, 
$\H = P_{\leq 2}(V) \cup \{ 123\}$ and $S = (V,\H)$. Then $S$ is not a TBRSC.
\ee

With respect to the equality $L(S^{\varepsilon}) = \varepsilon(S)$ in Theorem \ref{eqtr}, we can show it holds for all paving simplicial complexes:

\bp
\label{eqp}
Let $d \geq 0$ and $S = (V,\H) \in {\rm Pav}(d)$. Then $L(S^{\varepsilon}) = \varepsilon(S)$.
\ep

\proof
In the Appendix.
\qed

However, the next example, also analyzed in the Appendix, shows that equality may not hold.

\be
\label{ttnot}
Let $S = (V,\H)$ be the simplicial complex defined by $V = \{ 1,\ldots,5\}$ and
$$\H = (P_{\leq 2}(V) \setminus \{ 14,24,35\}) \cup \{ 123 \}.$$
Then $\varepsilon(S) \subset L(S^{\varepsilon})$.
\ee

\section{Low dimensions}

\bp
\label{tbone}
Every TBRSC of dimension 1 is boolean representable.
\ep

\proof
Let $S = (V,\H)$ be a TBRSC of dimension 1. For every $a \in V$, let $F_a = \{ a\} \cup \{ b \in V \setminus \{ a\} \mid ab \notin \H \}$. Every $a \in V$ is a transversal of the successive differences for the chain $\emptyset \subset V$ in $L(S)$. Suppose now that $a,b \in V$ are distinct and $ab \in \H$. Then $ab$  
is a transversal of the successive differences for the chain $\emptyset \subset F_a \subset V$, so it suffices to show that $F_a \in L(S)$. 

Since $S$ is a TBRSC, there exists a boolean matrix $M$ with column space $V$ such that, for every $X \in P_{\leq 2}(V)$, we have $X \in \H$ if and only if $X$ is $M$-independent. Since $P_1(V) \subseteq \H$, all the columns of $M$ are nonzero, so $X \in \H$ if and only if the columns of $X$ are different. Thus $F_a$ is the set of all $b \in V$ having columns in $M$ equal to the column of $a$. 

Let $X \in \H \cap 2^{F_a}$ and $p \in V \setminus F_a$. Then $|X| \leq 1$ and the column of $p$ is different, so $X \cup \{ p\} \in \H$ and so $F_a \in L(S)$ as required.
\qed 

Example \ref{truno} shows that Proposition \ref{tbone} fails for dimension 2, even in the paving case.

The next lemma features a class of matroids which is useful to build counterexamples.

\bl
\label{buildm}
Let $V$ be a finite nonempty set and let $\F \subseteq P_{2}(V) \cup P_3(V)$ be such that $F \cap F' = \emptyset$ for all distinct $F,F' \in \F$. Let $\H$ consist of all $X \in P_{\leq 3}(V)$ containing no element of $\F$. Then $(V,\H)$ is a matroid.
\el

\proof
In the Appendix.
\qed

The next result shows that, when it comes to separate BRSCs from TBRSCs, Example \ref{truno} has the minimum number of vertices.

\bp
\label{btbfive}
Every TBRSC with at most 5 vertices is boolean representable.
\ep

\proof
Let $S = (V,\H)$ be a TBRSC with $|V| \leq 5$.
In view of Proposition \ref{tbone}, we may assume that $\dim(S) \geq 2$.

Suppose first that $|V| = 3$. Then $S$ is the uniform matroid $U_{3,3}$, hence a BRSC.

Suppose next that $|V| = 4$. We may assume that $\dim(S) = 2$, otherwise $S = U_{4,4}$. If $S: = T_3(S')$ for some BRSC $S'$ of dimension 3, then $S = U_{3,4}$ and is therefore a BRSC.

Thus we may assume that $|V| = 5$. If $\dim(S) = 4$, then $S = U_{5,5}$ is a BRSC.
If $\dim(S) = 3$ and $S = T_4(S')$ for some BRSC $S'$ of dimension 4, then $S = U_{4,5}$ and is also a BRSC. Hence we may assume that $\dim(S) = 2$.

Suppose that $S$ is not a BRSC. Then $L(S) \subset \varepsilon(S)$. Let $Z \in \varepsilon(S) \setminus L(S)$. Comparing the definitions of $\varepsilon(S)$ and $L(S)$, we see that
 $\H \cap P_3(Z) \neq \emptyset$, hence we may take $a_1a_2a_3 \in \H \cap P_3(Z)$. Since $a_1a_2a_3 \in \H \subseteq \H^{\varepsilon}$, we may assume that there exists a chain $\emptyset = Z_0 \subset Z_1 \subset Z_2 \subset Z_3$ in $\varepsilon(S)$ such that $a_i \in Z_i \setminus Z_{i-1}$ for $i = 1,2,3$. By Lemma \ref{propt}(i), $a_1a_2a_3$ is also a transversal of the successive differences for the chain
$$\emptyset \subset Z_1 \cap Z \subset Z_2 \cap Z \subset Z_3 \cap Z$$ 
in $\varepsilon(S)$, hence there exists a chain 
\beq
\label{wol}
\emptyset = Z'_0 \subset Z'_1  \subset Z'_2 \subset Z'_3 \subset Z'_4 = V
\eeq 
in $\varepsilon(S)$. Since $|V| = 5$, there exists some $j \in 1234$ such that $|Z'_j \setminus Z'_{j-1}| = 2$. Without loss of generality, we may assume that 
$Z'_j = Z'_{j-1} \cup 12$. 

If $X \in P_3(V)$ does not contain $12$, then $X$ is a transversal of the successive differences for the chain (\ref{wol}), hence $X \in T_3(\H^{\varepsilon}) = \H$. Thus the only possible elements of $P_3(V) \setminus \H = P_3(V) \setminus \H^{\varepsilon}$ are $123, 124, 125$. 

If $12 \notin \H$, we have necessarily 
$$\H = \{ X \in P_{\leq 3}(V) \mid 12 \not\subseteq X \},$$
because any other 2-subset is contained in some element of $\H \cap P_3(V)$. By Lemma \ref{buildm}, $S$ is a matroid, hence boolean representable.

Thus we may assume that we have one of the following four cases:
\bi
\item[(C1)] $\H = P_{\leq 3}(V) \setminus \{ 123, 124, 125 \}$;
\item[(C2)] $\H = P_{\leq 3}(V) \setminus \{ 123, 124\}$;
\item[(C3)] $\H = P_{\leq 3}(V) \setminus \{ 123 \}$;
\item[(C4)] $\H = P_{\leq 3}(V)$.
\ei
Now (C3) and (C4) are clearly both matroids (hence BRSCs). We can show that (C1) is a BRSC by checking that $34,35,45$ are flats. Similarly, (C2) is a BRSC because $15, 34,35,45$ are flats. Therefore every TBRSC with 5 vertices is a BRSC.
\qed

\section{Join}

Given two simplicial complexes $S = (V,\H)$ and $S' = (V,\H')$ we define the {\em join} of $S$ and $S'$ as the simplicial complex
$$S \vee S' = (V,\H \cup \H').$$
Notice that given a simplicial complex $(V,\H)$, then $\H$ is just a down set of $2^V$ under inclusion. The down sets of $2^V$ form a lattice equal to the lattice of semigroup ideals of the monoid $(2^V, \cap)$, 
and this construction is precisely the join in this lattice.

\bp
\label{unfo}
Let $S = (V,\H)$ and $S' = (V,\H')$ be BRSCs with $|V| \leq 4$. Then $S  \vee S'$ is a BRSC.
\ep

\proof
If $\dim(S \vee S') = 1$, we may use Proposition \ref{tbone} and Theorem \ref{pavun}. The only other nontrivial case is $\dim(S  \vee S') = 2$. But it is easy to check \cite[Example 5.2.11]{RSm} that if $|V| = 4$ and $\dim(S) = 2$, then $S$ is a BRSC if and only if $|\H \cap P_3(V)| \neq 1$. It follows that if $S \vee S'$ is not a BRSC, then $S$ or $S'$ is not a BRSC.
\qed

The next example shows that neither BRSCs nor TBRSCs are closed under join when we consider 5 vertices (even at dimension $\leq 2$). We analyze this example in the Appendix.

\be
\label{nonun}
Let $V = 12345$. Let $S_1 = (V,P_{\leq 2}(V))$ and $S_2 = (V,\H_2)$ be defined by
$$\H_2 = \{ X \in P_{\leq 3}(V) \mid 12, 34 \not\subseteq X \}.$$ 
\ee

But things work out better in the paving case:

\bt
\label{pavun}
Let $d \geq 1$ and let $(V,\H), (V,\H') \in {\rm TBPav}(d)$. Then $(V,\H\cup \H') \in {\rm TBPav}(d)$.
\et

\proof
Let 
$$R = \{ Z \cap Z' \mid Z \in \varepsilon(\H),\, Z' \in \varepsilon(\H')  \}.$$
In view of Lemma \ref{propt}(i), $R$ is a Moore family. Hence $(V,{\rm Tr}(R))$ is a BRSC. We claim that 
\beq
\label{pavun1}
\H \cup \H' = T_{d+1}({\rm Tr}(R)).
\eeq

Let $X \in \H$. By Theorem \ref{eqtr}, there exists a chain
\beq
\label{pavun2}
Z_0 \subset Z_1 \subset \ldots \subset Z_n
\eeq
in $\varepsilon(\H)$ and an enumeration $x_1,\ldots,x_n$ of the elements of $X$ such that $x_i \in Z_i \setminus Z_{i-1}$ for every $i$. Since $V \in \varepsilon(\H')$, then (\ref{pavun2}) is also a chain in $R$, hence $X \in {\rm Tr}(R)$. But $|X| \leq d+1$, thus $\H \subseteq T_{d+1}({\rm Tr}(R))$ and also $\H' \subseteq T_{d+1}({\rm Tr}(R))$ by symmetry.

Conversely, let $X \in T_{d+1}({\rm Tr}(R))$. Since $\H,\H' \in \pav(d)$, we may assume that $|X| = d+1$. 
Then there exists a chain
\beq
\label{pavun3}
R_0 \subset R_1 \subset \ldots \subset R_{d+1}
\eeq
in ${\rm Tr}(R)$ and an enumeration $x_1,\ldots,x_{d+1}$ of the elements of $X$ such that $x_i \in R_i \setminus R_{i-1}$ for every $i$. 

Write $R_d = Z \cap Z'$ with $Z \in \varepsilon(\H)$ and $Z' \in \varepsilon(\H')$. Since $x_{d+1} \notin R_d$, we may assume that $x_{d+1} \notin Z$. Since $(V,\H) \in \tbpav(d)$ yields $P_{\leq d-1}(V) \subseteq \varepsilon(\H)$, then
$$\emptyset \subset x_1 \subset x_1x_2 \subset \ldots \subset x_1\ldots x_{d-1} \subset Z \subset V$$
is a chain in $\varepsilon(\H)$ having $X$ as a transversal of the successive differences. Thus $X \in \H$ by Theorem \ref{eqtr} and so (\ref{pavun}) holds. Note also that $P_{\leq d}(V) \subseteq \H \subseteq \H \cup \H'$. 

Therefore $(V,\H\cup \H') = T_{d+1}(V,{\rm Tr}(R)) \in {\rm TBPav}(d)$.
\qed

The next example, analyzed in the Appendix, shows that we cannot replace $\tbpav(d)$  by $\bpav(d)$ in Theorem \ref{pavun}.

\be
\label{ncu}
Let $V = 123456$, $\H = P_{\leq 2}(V) \cup \{ 123, 124, 125, 126\}$ and 
$$\H' = P_{\leq 2}(V) \cup 
\{ X \in P_3(V) \; {\big{\lvert}} \; |X \cap 46| = 1\}.$$
Then $(V,\H),(V,\H') \in {\rm BPav}(d)$ but $(V,\H\cup \H') \notin {\rm BPav}(d)$.
\ee

Let $V$ be a finite nonempty set and let $L \subseteq V$ be such that $2 \leq d \leq |L| < |V|$. We write
$$B_d(V,L) = (V,\B_d(V,L)) = (V,{\rm Tr}(P_{\leq d-1}(V) \cup \{ V,L\})).$$
This is easily seen to be equivalent to the following condition:
$$\B_d(V,L) = P_{\leq d}(V) \cup \{ X \in P_{d+1}(V) \; {\big{\lvert}} \; |X \cap L| = d\}.$$
If $V$ is clear from the context, we may omit $V$ from $\B_d(V,L)$ and $B_d(V,L)$.

\bl
\label{bvl}
Let $V$ be a finite nonempty set and let $L \subseteq V$ be such that $2 \leq d \leq |L| < |V|$. Then $B_d(V,L) \in {\rm BPav}(d)$.
\el

\proof
It is immediate that $P_{\leq d-1}(V) \cup \{ V,L\} \subseteq L(B_d(L))$, hence every $X \in \B_d(L)$ is a transversal of the successive differences for some chain in $L(B_d(L))$. Thus $B_d(L) \in \bpav(d)$.
\qed

We can now prove the following result, characterizing ${\rm TBPav}(d)$.

\bt
\label{ubvl}
Let $d \geq 2$ and $S = (V,\H) \in {\rm Pav}(d)$. Then the following conditions are equivalent:
\bi
\item[(i)] $S \in {\rm TBPav}(d)$;
\item[(ii)] $S = \vee \{ B_d(V,L) \mid L \in \L \}$ for some nonempty $\L \subseteq P_{\geq d}(V) \setminus \{ V\}$.
\ei
\et

\proof
(i) $\Rw$ (ii). Let $\L = (P_{\geq d}(V) \setminus \{ V\}) \cap \varepsilon(S)$. Since $\dim(S) = d$, we have $\L \neq \emptyset$.

Let $X \in \H$. Since $P_{\leq d}(V) \subseteq \B_d(L)$ for every $L \in \L$, we may assume that $|X| = d+1$. By Theorem \ref{eqtr}, there exists a chain
$$Z_0 \subset Z_1 \subset \ldots \subset Z_{d+1}$$
in $\varepsilon(S)$ and an enumeration $a_1,\ldots a_{d+1}$ of the elements of $X$ so that $a_i \in Z_i \setminus Z_{i-1}$ for $i = 1,\ldots,d+1$. Now $a_1\ldots a_{i} \in L(S) \subseteq \varepsilon(S)$ for $i = 0,\ldots,d-1$, hence $X$ is a transversal of the chain
$$\emptyset \subset a_1 \subset a_1a_2 \subset \ldots a_1\ldots a_{d-1} \subset Z_d\subset V$$
and so $X \in \B_d(T_d)$. Since $Z_d \in \L$, we get $\H \subseteq \vee \{ B_d(L) \mid L \in \L \}$.

The opposite inclusion is immediate.

(ii) $\Rw$ (i). By Lemma \ref{bvl} and Theorem \ref{pavun}.
\qed

In such a decomposition ($S = \vee \{ B_d(V,L) \mid L \in \L \}$), we may refer to the elements of $\L$ as {\em lines}.

The following lemma shows that the decomposition provided by Theorem \ref{ubvl} is not unique.

\bl
\label{bigshort}
Let $d \geq 2$ and let $V$ be a finite set with $|V| \geq d+1$. For every $a \in V$, we have
\beq
\label{bigshort1}
\B_d(V, V\setminus \{ a\}) =  \displaystyle\bigcup_{L \in \L_a} \B_d(V,L),
\eeq
where $\L_a = \{ L \in P_d(V) \mid a \in L\}$.
\el

\proof
It suffices to show that both sides of (\ref{bigshort1}) contain the same $X \in P_{d+1}(V)$. So let $X \in P_{d+1}(V)$.

Suppose that $X \in \B_d(V, V\setminus \{ a\})$. Then $|X \cap (V \setminus \{ a\})| = d$, hence $a \in X$. Take $b \in X \setminus \{ a \}$. Then $X \setminus \{ b \} \in \L_a$ and so $$X \in \B_d(V, X \setminus \{ b \}) \subseteq \displaystyle\bigcup_{L \in \L_a} \B_d(V,L).$$

Conversely, suppose that $X \in \B_d(V,L)$ with $L \in \L_a$. Since $|X| = d+1$ and $|L| = d$, we must have $X = L \cup \{ c \}$ for some $c \in V \setminus L$. Hence $a \in L \subset X$ yields $|X \cap (V\setminus \{ a\})| = d$ and $X \in \B_d(V, V\setminus \{ a\})$. Therefore (\ref{bigshort1}) holds as required.
\qed

We prove next a version of Theorem \ref{ubvl} for $\bpav(d)$.

\bt
\label{ubp}
Let $d \geq 2$ and $S = (V,\H) \in {\rm Pav}(d)$. Then the following conditions are equivalent:
\bi
\item[(i)] $S \in {\rm BPav}(d)$;
\item[(ii)] $S = \vee \{ B_d(V,L) \mid L \in L(S) \setminus (P_{\leq d-1}(V) \cup \{ V \}) \}$;
\item[(iii)] $S = \vee \{ B_d(V,L) \mid L \in \L \}$ for some nonempty $\L \subseteq P_{\geq d}(V) \setminus \{ V\}$ satisfying
\beq
\label{bdlbr1}
|L \cap L'| \leq d-1 \mbox{ for all distinct }L,L' \in \L.
\eeq
\ei
\et

\proof
(iii) $\Rw$ (i). 
Since $P_{\leq d}(V) \subseteq \H$, we have $P_{\leq d-1}(V) \subseteq L(S)$. Let $K \in \L$ and suppose that $X \in \H \cap 2^K$ and $p \in V \setminus K$. Since $P_{\leq d}(V) \subseteq \H \subseteq P_{\leq d+1}(V)$, we may assume that $|X| = d$ or $d+1$. 

Suppose that $|X| = d+1$. Since $X \in \H = \displaystyle\cup_{L \in \L} \B_d(L)$, we have $X \in \B_d(L)$ for some $L \in L$. Thus $|X \cap L| = d$ and so $|K \cap L| \geq d$. In view of (\ref{bdlbr1}), we get $K = L$, hence $X \subseteq L$, a contradiction since $|X| = d+1$ and $|X \cap L| = d$. Therefore
$|X| = d$, hence $X \cup \{ p \} \in \B_d(K) \subseteq H$ and so $K \in L(S)$. 

Let $a_1, \ldots, a_{d-1} \in V$ be distinct. Then
\beq
\label{bdlbr2}
\emptyset \subset a_1 \subset a_1a_2 \subset \ldots \subset a_1\ldots a_{d-1} \subset V
\eeq
is a chain in $L(S)$. If $a_1, \ldots, a_{d-1} \in L \in \L$, then (\ref{bdlbr2}) can be refined to
\beq
\label{bdlbr3}
\emptyset \subset a_1 \subset a_1a_2 \subset \ldots \subset a_1\ldots a_{d-1} \subset L \subset V,
\eeq
also a chain in $L(S)$. It is easy to check that every $X \in \H$ is a partial transversal of the successive differences for some chain of type (\ref{bdlbr2}) or (\ref{bdlbr3}), hence 
$S$ is boolean representable.

(i) $\Rw$ (ii). 
Since $P_{\leq d-1}(V) \cup \{ V \} \subseteq L(S)$ and the maximum length of a chain in $L(S)$ is $d+1$, it follows easily that the maximal chains in $L(S)$ must be of the form (\ref{bdlbr2}) or (\ref{bdlbr3}), with $L \in L(S) \setminus (P_{\leq d-1}(V) \cup \{ V \})$. Thus (ii) holds.

(ii) $\Rw$ (iii). 
Suppose that $L,L' \in L(S)$ are distinct and satisfy $|L \cap L'| \geq d$. We may assume that $L\cap L' \subset L$. Let $a_1,\ldots,a_{d-1} \in L$ be distinct. Since $P_{\leq d-1}(V) \subseteq L(S)$, we get a chain of length $d+2$
$$\emptyset \subset a_1 \subset a_1a_2 \subset \ldots \subset a_1\ldots a_{d-1} \subset L\cap L' \subset L \subset V$$
in $L(S)$, contradicting $\dim\H\, = d$. Therefore $|L \cap L'| \leq d-1$.
\qed

\bc
\label{smallpav}
Let $d \geq 2$ and and let $V$ be a finite set with $|V| \geq d+1$. Let $\emptyset \neq \L \subseteq 2^V$ be such that $|L| \in \{ d,d+1, |V|-1\}$ for every $L \in L$. Then $\vee\{ B_d(V,L) \mid L \in \L \}$ is boolean representable.
\ec

\proof
In view of Lemma \ref{bigshort}, we may assume that $|L| \in \{ d,d+1 \}$ for every $L \in \L$. Let $L,L' \in \L$ be distinct.
It is easy to check that
\beq
\label{smallpav1}
\B_d(L) \cup \{ L\} = \cup\{ \B_d(L \setminus \{ a \}) \mid a \in L \}
\eeq
holds for every $L \in P_{d+1}(V)$. 
Now we may use (\ref{smallpav1}) for replacing $\L$ by some equivalent $\L' $ satisfying (\ref{bdlbr1}):
\bi
\item[(1)] 
If $L,L' \in \L \cap P_{d+1}(V)$ are such that $|L \cap L'| = d$, we replace $\B_d(L) \cup \B_d(L')$ by
$$(\cup\{ \B_d(L \setminus \{ a \}) \mid a \in L \}) \cup (\cup\{ \B_d(L' \setminus \{ a \}) \mid a \in L' \}).$$
\item[(2)] 
If $L \in \L \cap P_{d+1}(V)$ and $L' \in P_d(L)$, we replace $\B_d(L)\cup \B_d(L')$ by
$\cup\{ \B_d(L \setminus \{ a \}) \mid a \in L \}$.
\ei
Indeed, these replacements are legitimate in view of (\ref{smallpav1}), and each such replacement decreases the number of $L \in \L \cap P_{d+1}(V)$. Eventually, we end up with some $\L' $ satisfing (\ref{bdlbr1}).
By Theorem \ref{ubp}, our complex is boolean representable.
\qed

Another way of ensuring closure under join is by restricting the type of complexes in $\bpav(d)$.
We define, for every $d \geq 2$ and every finite set $V$ with at least $d+2$ elements,
$$\Y(V) = \{ (V,\H) \in \bpav(d) \mid (V,\H)\mbox{
 has no restriction isomorphic to }U_{d,d+2} \}.$$
 
\bp
\label{uniony}
Let $d \geq 2$ and let $V$ be a finite set with at least $d+2$ elements. Then
$$(V,\H_1), (V,\H_2) \in \Y(V) \mbox{ implies } (V,\H_1 \cup \H_2) \in \Y(V).$$
\ep

\proof
We show that
\beq
\label{uniony2}
\flatx(V,\H_i) \subseteq P_{\leq d+1}(V) \cup \{ V \}.
\eeq
Let $i \in \{ 1,2\}$ and $F \in L(V,\H_i)$. Suppose that $d+1 < |F| < |V|$. Since $P_d(V) \subseteq \H_i$ and no restriction of $(V,\H_i)$ to a $d$-subset of $F$ is isomorphic to $U_{d,d+2}$, there exists some $X \in P_{d+1}(F) \cap \H_i$. But then $F$ contains a facet of $(V,\H_i)$ and so $F = V$, a contradiction. Therefore (\ref{uniony2}) holds.

Now let 
$$\L_i = \{ F \in L(V,\H_i) \mid d \leq |F| < |V| \}.$$
By (\ref{uniony2}), we have
\beq
\label{uniony1}
\L_i = \{ F \in L(V,\H_i) \; {\big{\lvert}} \; |F| \in \{ d,d+1\} \, \}.
\eeq
Since $(V,\H_i) \in \bpav(d)$, we have $\H_i = \displaystyle\bigcup_{L \in \L_i} \B_d(V,L)$ by Theorem \ref{ubp}. Thus $\H_1 \cup \H_2 = \cup\{ \B_d(V,L) \mid L \in \L_1 \cup \L_2 \}$ and so $(V,\H_1 \cup \H_2) \in \bpav(d)$ by (\ref{uniony1}) and Proposition \ref{smallpav}. 

Suppose that there exists some $W \in P_{d+2}(V)$ such that $P_{d+1}(W) \cap (\H_1 \cup \H_2) = \emptyset$. Then $P_{d+1}(W) \cap \H_1 = \emptyset$, contradicting $(V,\H_1) \in \Y(V)$. Therefore $(V,\H_1 \cup \H_2) \in \Y(V)$. 
\qed

Finally, we prove that, when we start with a complex $S \in {\rm Pav}(d)$, there is a largest paving TBRSC contained in $S$, and it is precisely $T_{d+1}(S^{\varepsilon})$.

\bt
\label{sczero}
Let $S = (V,\H) \in {\rm Pav}(d)$ and 
$$\tau(S) = \{ (V,\H') \in {\rm TBPav}(d) \mid \H' \subseteq \H\} \cup \{ (V,P_{\leq d}(V)) \}.$$
Then:
\bi
\item[(i)]
there exists some (unique) $S^0 = (V,\H^0) \in \tau(S)$ such that $\H' \subseteq \H^0$ for every $(V,\H') \in \tau(S)$;
\item[(ii)]
$S^0 = T_{d+1}(S^{\varepsilon})$.
\ei
\et

\proof
(i) Let $\H^0 = \cup\{ \H' \mid (V,\H') \in \tau(S) \}$. Clearly, $P_{\leq d}(V) \subseteq \H'$ for every $(V,\H') \in \tau(S)$. In view of Theorem \ref{pavun}, it follows that $(V,\H^0) \in \tau(S)$ and we are done.

(ii) By Lemma \ref{inc}, $T_{d+1}(S^{\varepsilon})$ is a TBRSC and $T_{d+1}(\H^{\varepsilon}) \subseteq \H$. Since $S \in \pav(d)$, we have $P_{\leq d-1}(V) \cup \{ V\} \subseteq \varepsilon(S)$ and so $P_{\leq d}(V) \subseteq \H^{\varepsilon}$. Hence $T_{d+1}(S^{\varepsilon})$ is paving and $T_{d+1}(S^{\varepsilon}) \in \tau(S)$. Therefore $T_{d+1}(\H^{\varepsilon}) \subseteq \H^0$.

To prove the opposite inclusion, we may assume that $\dim(S^0) = d$ (otherwise $\H^0 = P_{\leq d}(V) \subseteq \H^{\varepsilon}$ and we are done). It follows from Theorem \ref{eqtr} that $S^0 = T_{d+1}((S^0)^{\varepsilon})$, so it suffices to show that 
$(\H^0)^{\varepsilon} \subseteq \H^{\varepsilon}$, which follows itself from $\varepsilon(S^0) \subseteq \varepsilon(S)$. We prove the latter inclusion.

Let $Z \in \varepsilon(S^0)$. Then
$$\forall X \in \H^0 \cap P_{\leq d}(Z)\; \forall p \in V
\setminus Z \hspace{.3cm} X \cup \{ p \} \in \H^0.$$
Suppose that $X \in \H \cap P_{\leq d}(Z)$ and $p \in V
\setminus Z$. Since $S^0 \in \pav(d)$, we have $X \in \H^0 \cap P_{\leq d}(Z)$ and $Z \in \varepsilon(S^0)$ yields $X \cup \{ p \} \in \H^0 \subseteq \H$. Thus $Z \in \varepsilon(S)$ and so $\varepsilon(S^0) \subseteq \varepsilon(S)$ as required.
\qed

The next example, analyzed in the Appendix, shows that the paving requirement for subcomplexes cannot be removed from the definition of $\tau(S)$.

\be
\label{scznp}
Let $S = (V,\H) \in {\rm Pav}(3)$ be defined by
$V = \on{7}$ and $\H = P_{\leq 3}(V) \cup \{ 1abc \mid a \in 23, b \in 45, c \in 67 \}$. Then $S$ has no largest truncated boolean representable subcomplex.
\ee

\section{On $\tbpav(d) \setminus \bpav(d)$}

We proved in Proposition \ref{btbfive} that we need at least 6 points to separate $\tbpav(2)$ from $\bpav(2)$. This section starts with a full account of the 6 point case.

\bp
\label{six}
Up to isomorphism, the complexes with 6 points in $\tbpav(2) \setminus \bpav(2)$ are of the form $(123456,\H)$ for:
\bi
\item[(1)] $\H = \B_2(1234) \cup \B_2(12)$;
\item[(2)] $\H = \B_2(1234) \cup \B_2(12) \cup \B_2(15)$;
\item[(3)] $\H = \B_2(1234) \cup \B_2(12) \cup \B_2(15) \cup \B_2(25)$;
\item[(4)] $\H = \B_2(1234) \cup \B_2(12) \cup \B_2(35)$;
\item[(5)] $\H = \B_2(1234) \cup \B_2(12) \cup \B_2(15) \cup \B_2(35)$. 
\ei
Moreover, all the above 5 cases are nonisomorphic.
\ep

\proof
In the Appendix.
\qed

Ordering the sets of faces through inclusion, we can build the following diagram
$$
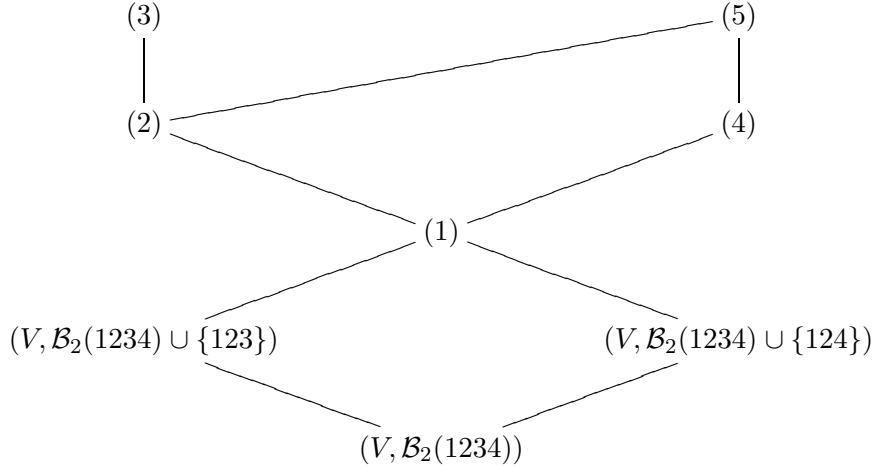
\xymatrix{
(3) \ar@{-}[d]&&(5) \ar@{-}[d] \ar@{-}[dll]\\
(2)\ar@{-}[dr]&&(4) \ar@{-}[dl]\\
&(1) \ar@{-}[dr] \ar@{-}[dl] &\\
(V,\B_2(1234) \cup \{123\}) \ar@{-}[dr]&&(V,\B_2(1234) \cup \{124\}) \ar@{-}[dl]\\
&(V,\B_2(1234))&
}$$
The missing triangles in the three lowest elements are respectively
$${124, 134, 234, 156, 256, 356, 456},\hspace{3cm}{123, 134, 234, 156, 256, 356, 456}$$
$${123, 124, 134, 234, 156, 256, 356, 456}$$
hence all the edges correspond to covering relations (check the enumeration of the missing triangles for (1)--(5) in the proof of Proposition \ref{six}).

We note that:
\bi
\item
$(V,\B_2(1234)) \in {\rm BPav}(2)$ by Lemma \ref{bvl}.
\item
$(V,\B_2(1234) \cup \{123\}) \notin {\rm TBPav}{2}$. Indeed, suppose that there exist $Z \in \varepsilon(\B_2(1234) \cup \{123\})$ such that $|Z \cap 123| = 2$. Since $124, 134, 234 \notin \H$, we successively get $4 \in Z$ and $1234 \subseteq Z$, a contradiction. In view of Theorem \ref{eqtr}, this implies $(V,\B_2(1234) \cup \{123\}) \notin {\rm TBPav}{2}$. 
\item
$(V,\B_2(1234) \cup \{124\}) \notin {\rm TBPav}{2}$. Similar to the preceding case.
\item
No simplicial complex isomorphic to (4) embeds in (3). To prove this, recall the missing triangles in (3) and (4). We can check that $3x$ is contained in a missing triangle of (3) for every $x \neq 3$. Similarly, $4y$ is contained in a missing triangle of (3) for every $y \neq 4$. Suppose that $\varphi \in S_6$ is such that the isomorphic image of (3) through $\varphi$ (call it (3'')) has (4) as subcomplex. Then the missing triangles of (3'') are a proper subset of the missing triangles of (4). Hence $(3\p)x$ is contained in a missing triangle of (4) for every $x \neq 3\p$, and $(4\p)y$ is contained in a missing triangle of (4) for every $y \neq 4\p$. However, only $4$ satisfies this property, yielding $3\p = 4 = 4\p$, a contradiction.
\ei

Note also that an arbitrary $S \in {\rm TBPav}(d) \setminus {\rm BPav}(d)$ needs not having a restriction isomorphic to $U_{d,d+2}$. The complexes featuring Proposition \ref{six} constitute all counterexamples for $d = 2$.

We intend now to show that ${\rm TBPav}(d) \setminus {\rm BPav}(d)$ is in some sense finitely generated. We start with a couple of lemmas.

Let $\tbr$ (respectively $\tbp$) denote the class of all finite truncated boolean representable simplicial complexes (respectively finite paving truncated boolean representable simplicial complexes).

A class of simplicial complexes closed under isomorphism and restriction is called a {\em prevariety}. For details on prevarieties, see \cite[Sections 8.4 and 8.5]{RSm}.

\bl
\label{pret}
The classes $\tbr$ and $\tbp$ are prevarieties of simplicial complexes.
\el

\proof
In the Appendix.
\qed

Let $S \in \tbpav(d) \setminus \bpav(d)$. By Lemma \ref{pret}, every restriction of $S$ is in $\tbp$ (with possibly lower dimension). We say that $S$ is
 {\em minimal} if every proper restriction of $S$ is boolean representable.
 
\bl
\label{comi}
Let $d \geq 2$. Then the maximum number of vertices for a minimal $S \in {\rm TBPav}(d) \setminus {\rm BPav}(d)$ is $(d+1)(d+2)$.
\el

\proof
Let $S = (V,\H) \in \tbpav(d) \setminus \bpav(d)$ be minimal. Hence $S \notin \bpav(d)$ but every proper restriction of $S$ is boolean representable. By \cite[Theorem 8.5.2(ii)]{RSm}, we get $|V| \leq (d+1)(d+2)$. 

Now we consider the {\em Swirl}, the simplicial complex defined in the proof of \cite[Theorem 8.5.2(ii)]{RSm}, where it is proved that every proper restriction of this complex is boolean representable, but the Swirl is not. The Swirl is defined as follows:

Let $A = \{ a_0,\ldots, a_d\}$ and $B_i =  \{ b_{i0},\ldots, b_{id}\}$
for $i = 0,\ldots,d$. Write also $A_i = A \setminus \{ a_i
\}$ and 
$$C_i = P_{d+1}(A_i \cup (B_i\setminus \{ b_{i0} \})) \cup \{ B_i \}.$$
We define $$V = A \cup \displaystyle\bigcup_{i=0}^d B_i,\quad
\H = P_{\leq d+1}(V) \setminus \displaystyle\bigcup_{i=0}^d C_i.$$
It is easy to check that all the $X \in \H \cap P_{d+1}(V)$ fall into four cases (not necessarily disjoint):
\bi
\item[(a)]  
there exist $b_{ij}, b_{k\ell} \in X$ with $i \neq k$;
\item[(b)] 
there exist $b_{i0}, a_j \in X$;
\item[(c)] 
there exist $b_{ij}, a_i \in X$ with $j > 0$;
\item[(d)] 
$X = a_0\ldots a_d$.
\ei
Define
$$\begin{array}{lll}
\L&=&\{ L \in P_d(V) \mid \mbox{ there exist some $b_{ij}, b_{k\ell} \in L$ with $i \neq k$} \}\\
&\cup&\{ L \in P_d(V) \mid \mbox{ there exist some $b_{i0}, a_j \in L$} \}\\
&\cup&\{ L \in P_d(V) \mid \mbox{ there exist some $b_{ij}, a_i \in L$ with $j > 0$} \}\\
&\cup&\{ A_i \cup B_i \mid i = 0,\ldots,d\}.
\end{array}$$
It is straightforward to check that $\H = \cup\{ \B_d(L) \mid L \in \L \}$, hence $S \in \tbpav(d)$ by Theorem \ref{ubvl}. Since
$|V| = (d+1)(d+2))$, we have found some minimal $S \in \tbpav(d) \setminus \bpav(d)$ with $(d+1)(d+2)$ vertices as required.
\qed

Let $\V$ be a prevariety of simplicial complexes. We say that $\V$ is {\em finitely based} if there exists some $m \geq 1$ such that every simplicial complex not in $\V$ admits a restriction not in $\V$ with at most $m$ vertices.

Given a prevariety $\V$ of simplicial complexes and $d \in
\mathbb{N}$, we define the prevariety
$$\V_d = \{ S \in \V \mid \dim(S) \leq d \}.$$

Let $\cal{BP}$ denote the class of all finite paving boolean representable simplicial complexes. By \cite[Theorem 8.5.2]{RSm}, ${\cal{BP}}_d$ is finitely based for every $d \geq 1$. Since ${\cal{TBP}}_1 = {\cal{BP}}_1$ by Proposition \ref{tbone}, it follows that ${\cal{TBP}}_1$ is finitely based.

\bt
\label{nfb}
${\cal{TBP}}_2$ is not finitely based.
\et

\proof
It suffices to build arbitrary large simplicial complexes not in ${\cal{TBP}}_2$ with all proper restrictions in ${\cal{TBP}}_2$.

Suppose that $S = (V,\H) \in \pav(2)$. Then $P_{\leq 1}(V) \subseteq L(S) \subseteq \varepsilon(S)$, and so
\begin{eqnarray}\nonumber
\label{nfb1}
S \in \tbp_2 \mbox{ if and only if, for every $X \in \H \cap P_3(V)$,}\hspace{3cm}\\
\hspace{2cm}\mbox{ there exists some $Z \in \varepsilon(S)$ such that }|X \cap Z| = 2.
\end{eqnarray}

Let $n \geq 6$ and take as vertex set
$$V = \{ x_0, \ldots, x_n, y_0, \ldots, y_6, z_0, \ldots, z_6 \},$$
where we identify
$$x_0 = y_0 = z_6, \quad
x_1 = z_0 = y_6, \quad 
y_1 = z_1 = x_n.$$
Let
$$T = \{ x_ix_{i+1}x_{i+2} \mid i = 0,\ldots, n-2\} \cup 
\{ y_iy_{i+1}y_{i+2} \mid 0 \leq i \leq 4 \} \cup \{ z_iz_{i+1}z_{i+2} \mid 0 \leq i \leq 4 \},$$
$\H = P_{\leq 3}(V) \setminus T$ 
and $S = (V,\H)$. We show next that:
\beq
\label{nfb2}
\mbox{for every $X \in (\H \cap P_3(V)) \setminus \{ x_0x_1y_1\}$, there exists some $F \in L(S)$ such that }|X \cap F| = 2.
\eeq
Indeed, such an $X$ contains necessarily some element of $V \setminus x_0x_1y_1$. Without loss of generality, we may assume that this element is among $x_2,\ldots,x_{n-1}$ (the other cases follow by symmetry). 

Suppose that $X \subset x_0\ldots x_n$, say $X = x_ix_jx_k$ with $i < j < k$. Since $X \in \H$, then $i,j,k$ are not consecutive integers. If $k < n$, then $k > 1$ and $k-i > 2$, hence $x_ix_k \in L(S)$ and we are done. Thus we may assume that $k = n$. If $i > 1$, then $k-i > 2$, hence $x_ix_k \in L(S)$ and we are done. Thus we may assume that $i \leq 1$. Since $k = n$, this implies $2 \leq j \leq n-1$. Since $n \geq 6$, we get either $k-j > 2$ (yielding $x_jx_k \in L(S)$) or $j-i > 2$ (yielding $x_ix_j \in L(S)$). 

Hence we may assume that at least one of the other elements of $X$ (say $a$) is not of the form $x_j$. Let $i \in \{ 2,\ldots,n-1\}$ be such that $x_i \in X$. It is easy to check that $x_ia \in L(S)$. Therefore (\ref{nfb2}) holds.

Next we show that 
\beq
\label{nfb3}
\mbox{for every $Z \in \varepsilon(S)$, $|Z \cap \{ x_0x_1y_1\}| \neq 2$.}
\eeq

Let $Z \in \varepsilon(S)$ and assume that $|Z \cap \{ x_0x_1y_1\}| \geq 2$. Assume first that $x_0,x_1 \in Z$. Since $x_ix_{i+1}x_{i+2} \notin H$ for $i = 0, \ldots, n-2$, we get successively $x_2 \in Z,\ldots, x_n \in Z$. Since $x_n = y_1$, we get $x_0x_1y_1 \subseteq Z$.

Suppose now that $x_0,y_1 \in Z$. Since $x_0 = y_0$, we use the same argument to deduce that $y_2 \in Z,\ldots, y_6 \in Z$. Since $y_6 = x_1$, we get $x_0x_1y_1 \subseteq Z$.

Finally, suppose that $x_1,y_1 \in Z$. Since $x_1 = z_0$ and $y_1 = z_1$, we use the same argument to deduce that $z_2 \in Z,\ldots, z_6 \in Z$. Since $z_6 = x_0$, we get $x_0x_1y_1 \subseteq Z$ and (\ref{nfb3}) is proved.

In view of (\ref{nfb1}), it follows from (\ref{nfb3}) that $S \notin \tbp_2$.

Fix now $v \in V$ and write $W = V\setminus \{ v\}$. We must show that $S|_{W} \in \tbp_2$ (since $\tbp_2$ is closed under restrictions, this implies that $S|_{W'} \in \tbp_2$ for any $W' \subset V$). 

Since $S \in \pav(2)$, we only need to show that the righthand side of (\ref{nfb1}) holds when we replace $S$ by $S|_W$. Let $X \in \H \cap P_3(W)$. Suppose first that $X \neq x_0x_1y_1$. By (\ref{nfb2}),
there exists some $F \in L(S)$ such that $|X \cap F| = 2$. It follows that $F \cap W \in L(S|_W)$.  Since $|X \cap (F \cap W)| = 2$, the desired condition is satisfied if $X \neq x_0x_1y_1$. 

Thus we may assume that $X = x_0x_1y_1$. It follows that either $v = x_i$ with $2 \leq i \leq n-1$ or $v = y_j$ or $v = z_j$ with $2 \leq j \leq 5$.

Suppose that $v = x_i$. Let $Z = x_0\ldots x_{i-1}$. It is immediate that $Z \in \varepsilon(S|_W)$ and $|Z \cap x_0x_1y_1| = 2$. If $v = y_j$ (respectively $v = z_j$), we take $Z = y_0\ldots y_{j-1}$ (respectively $Z = z_0\ldots z_{j-1}$). Therefore, in view of (\ref{nfb1}), we get $S|_{W} \in \tbp_2$ as required.
\qed

\section{The pure core}

Let $S = (V,\H)$ be a simplicial complex of dimension $d$. 
We say that $S$ is {\em pure} if all the facets of $S$ have dimension $d$.

We define $\pure(S) = (V',\H')$ by
$$V' = \cup(\H \cap P_{d+1}(V)),\quad \H' = \cup\{ 2^X \mid X \in \H \cap P_{d+1}(V) \}.$$
It is immediate that $\pure(S)$ is the largest pure subcomplex of $S$, also called the {\em pure core} of $S$.

This section is devoted to the following questions, all related to the concept of pure core:

\bq
\label{pc}
Let $S$ be a BRSC. Is {\rm pure}$(S)$ a BRSC?
\eq

\bq
\label{tpc}
Let $S$ be a BRSC and let $k \geq 1$. Is {\rm pure}$(T_k(S))$ a BRSC?
\eq

\bq
\label{tpct}
Let $\H$ be a BRSC and let $k \geq 1$. Is {\rm pure}$(T_k(S))$ a TBRSC?
\eq

Note that Problem \ref{tpct} admits the equivalent statement:
\beq
\label{tpct1}
\mbox{If $S$ is a TBRSC and $k \geq 1$, is {\rm pure}$(T_k(S))$ a TBRSC?}
\eeq
For the nontrivial implication, let $S$ be a TBRSC of rank $r$ and assume that Problem \ref{tpct} has a positive answer. Since $S$ is a TBRSC, we have $S = T_r(S')$ for some BRSC $S'$. 

Suppose first that $k \geq r$. Then $T_k(S) = S$, hence we must show that  {\rm pure}$(S)$ is a TBRSC. Since $S = T_r(S')$, our goal follows from the answer of Problem \ref{tpct} for $S'$ and $r$.

Assume now that $k < r$. It is easy to check that $T_k(S) = T_k(S')$. Since the answer of Problem \ref{tpct} implies that {\rm pure}$(T_k(S'))$ is a TBRSC, then {\rm pure}$(T_k(S))$ is a TBRSC and so (\ref{tpct1}) has also a positive answer.

Since matroids are pure and closed under restriction \cite{Oxl}, all questions have positive answers for matroids. We show that none of them admits a positive answer in general, but we establish particular cases.

The following example, analyzed in the Appendix, answers Problem \ref{pc} negatively for dimension 3. It also 
answers Problems \ref{tpc} and \ref{tpct} for dimension 3 and $k = 4$.

\be
\label{pobe}
Let $V = \{ 1,\ldots,7\}$ and  $R =  \{ \emptyset, 1,3, 5,12, 56, 356, 1234, V\}$. Then $S = (V,{\rm Tr}(R))$ is a BRSC but ${\rm pure}(S)$ is not a TBRSC.
\ee

However, we can find positive answers for all the problems in particular cases as we shall see. The following lemma will prove useful:

\bl
\label{wool}
Let $S = (V,\H)$ be a simplicial complex and let $I,J \in \H$ be such that $I \subseteq \oo{J}$. Then there exists some $I' \in \H$ such that $I \subseteq I'$ and $\oo{I'} = \oo{J}$.
\el

\proof
In the Appendix.
\qed

Matroids admit a wide variety of characterizations. One of them is the following: a simplicial complex $M = (V,\H)$ is a matroid if and only if
\beq
\label{mrf}
\mbox{for all $X,Y \in \H$, $\oo{X} = \oo{Y}$ implies $|X| = |Y|$.}
\eeq
Indeed, if $\oo{X} = \oo{Y}$ and $|X| > |Y|$, it follows easily from the exchange property that $Y \cup \{ x \} \in \H$ for some $x \in X \setminus Y$, hence $\oo{Y \cup \{ x \}} = \oo{Y}$. 
In  the case of matroids, the enumeration in (\ref{derby3}) can be chosen arbitrarily (see \cite{Oxl}). By taking $x$ as last, we reach a contradiction. Thus $|X| = |Y|$.

Conversely, suppose that (\ref{mrf}) holds. Let $I,J \in \H$ be such that $|I| = |J| +1$. Suppose that $J\cup \{ i \} \notin \H$ for every $i \in I \setminus J$. Then $I \subseteq \oo{J}$ and so by Lemma \ref{wool} this contradicts (\ref{mrf}). Thus $M$ satisfies the exchange property and is therefore a matroid.

We define a simplicial complex $S = (V,\H)$ to be a {\em near-matroid} if 
$$\oo{X} = \oo{Y} \subset V \mbox{ implies }|X| = |Y|$$
for all $X,Y \in \H$. The {\em rank function} $\rho: L(S) \setminus \{ V \} \to \mathbb{N}$ is defined by
$$F\rho = |X|, \mbox{ where $X \in \H$ is such that }\oo{X} = F.$$
Note that  such an $X$ exists by \cite[Proposition 4.2.4]{RSm}.

It follows from (\ref{mrf}) that every matroid is a near-matroid. The following result shows that the converse fails, in fact a near matroid needs not be boolean representable.

\bp
\label{bpnm}
Let $S$ be a simplicial complex of dimension $d \geq 0$.
\bi
\item[(i)] If $S$ is paving, then $S$ is a near-matroid.
\item[(ii)] If $S$ is boolean representable and $d \leq 2$, then $S$ is a near-matroid.
\ei
\ep

\proof
(i) Write $S = (V,\H)$ and suppose that $X,Y \in \H$ are such that $\oo{X} = \oo{Y} \subset V$. By 
\cite[Proposition 4.2.3]{RSm}, $X$ and $Y$ are not facets (it is easy to check that the closure of a facet must be $V$). Suppose that $|X| < d$. Since $P_d(V) \subseteq \H$, it follows that $\oo{X} = X$, so in this case we get indeed $Y = X$. Thus we may assume by symmetry that $|X|,|Y| \geq d$. Since $X$ and $Y$ are not facets, then $|X| = d = |Y|$ and so $S$ is a near-matroid.

(ii) Let $X,Y \in \H$ be such that $\oo{X} = \oo{Y} \subset V$. Since $\oo{\emptyset} = \emptyset$ and the closure of a facet is $V$, we may assume that $X,Y \notin \fct(S) \cup \{ \emptyset \}$.

Assume first that $X = \{ a \}$. Let $M$ be an $R \times V$ boolean matrix representing $S$. If $C_a = \{ b \in V \mid M[R,a] = M[R,b] \}$, then $C_a \subseteq \oo{a}$. Moreover, $J \in \H \cap 2^{C_a}$ implies $|J| \leq 1$, and since all the columns of $M$ are nonzero (we have $P_1(V) \subseteq \H$), it follows that $\oo{a} = C_a$. Thus $\oo{Y} = \oo{X} = C_a$ yields $|Y| = 1 = |X|$.

Assume now that $|X| > 1$. By the previous case, 
 we also have $|Y| > 1$. Since $\dim(S) \leq 2$ and $X,Y \notin \fct(S)$, we have necessarily $|X| = 2 = |Y|$. Therefore $S$ is a near-matroid.
\qed

Before discussing boolean representable near-matroids, we present two lemmas.

\bl
\label{opre}
Let $S = (V,\H)$ be a 
near-matroid and let $F,F' \in L(S)$ be such that $F \subset F' \subset V$. Then $F\rho < F'\rho$.
\el

\proof
In the Appendix.
\qed

\bl
\label{step}
Let $S = (V,\H)$ be a 
near-matroid and let $F,F' \in L(S)$ be such that $F \subset F' \subset V$. Let $a_1 \in F'\setminus F$ and $k = F'\rho -F\rho$. Then there exist $a_2,\ldots a_k \in V$ such that 
$$F \subset \oo{F\cup a_1} \subset \oo{F\cup a_1a_2} \subset \ldots \subset \oo{F\cup a_1\ldots a_k} = F'.$$
\el

\proof
In the Appendix.
\qed  

\bt
\label{brnm}
Let $S$ be a boolean representable near-matroid and let $k \geq 0$. Then:
\bi
\item[(i)] $T_k(S)$ is a BRSC;
\item[(ii)] ${\rm pure}(T_k(S))$ is a BRSC.
\ei
\et

\proof
(i) Write $S = (V,\H)$ and let $\oo{X}$ denote the closure of $X \subseteq V$ in $L(S)$. We define
$$\F_k = \{ F \in L(S) \mid F\rho < k\} \cup \{ V\}.$$
Since $L(S)$ is closed under intersection, it follows from Lemma \ref{opre} that $\F_k$ is a Moore family. Hence $(V,{\rm Tr}(\F_k))$ is a BRSC. We show that $T_k(S) = (V,{\rm Tr}(\F_k))$.

Let $X \in T_k(\H)$ and let $s = |X|$. Then there exists an enumeration $a_1,\ldots,a_s$ of the elements of $X$ such that 
$$\oo{a_1} \subset \oo{a_1a_2} \subset \ldots \subset \oo{a_1\ldots a_s}.$$
Hence $X$ is a transversal of the successive differences for
$$\emptyset \subset \oo{a_1} \subset \oo{a_1a_2} \subset \ldots \subset \oo{a_1\ldots a_{s-1}} \subset V,$$
which is a chain in $\F_k$. Thus $X \in {\rm Tr}(\F_k)$.

Conversely, assume that $X \in {\rm Tr}(\F_k)$. Since $\F_k \subseteq L(S)$, it follows that $X \in \H$.

Suppose that $|X| > k$. Since $X \in {\rm Tr}(\F_k)$, there exist some $F \in \F_k$ and $x \in X$ such that $F \cap X = X \setminus \{ x\}$. But $F = \oo{Y}$ for some $Y \in T_{k-1}(\H)$. Hence $X \setminus \{ x\} \subseteq \oo{Y}$ and by Lemma \ref{wool} there exists some $Z \in \H$ such that $X \setminus \{ x\} \subseteq Z$ and $\oo{Z} = \oo{Y} = F$. But then $|Z| \geq |X \setminus \{ x\}| \geq k > |Y|$, a contradiction since $S$ is a near-matroid. Thus $|X| \leq k$ and so $T_k(\H) = {\rm Tr}(\F_k)$ as claimed.

(ii) Let $\F'_k$ denote the set of all flats of $S$ occurring in chains of the form
$$F_0 \subset F_1 \subset \ldots \subset F_k$$
in $\F_k$. We claim that $\F'_k$ is a Moore family.

Let $F,F' \in \F'_k$. We may assume that $F,F' \neq V$. We have $F \cap F' \in \F_k$ since $\F_k$ is a Moore family. Since $F \in \F'_k \setminus \{ V \}$, there exists some $F'' \in L(S)$ such that $F \subseteq F''$ and $F''\rho = k-1$. Now we apply Lemma \ref{step} to both inclusions $\emptyset \subseteq F\cap F' \subseteq F''$. This ensures that $F\cap F'$ will appear in some chain of flats of length $k$ in $L(S)$ of the form
$$\emptyset \subset \ldots \subset F'' \subset V.$$
Since $F''\rho = k-1$, it follows from Lemma \ref{opre} that this is in fact a chain in $\F_k$ and therefore in $\F'_k$. Thus $F\cap F' \in \F'_k$. Since $V \in \F'_k$, then $\F'_k$ is a Moore family. Writing $V' = \cup \F'_k$, it follows that
$(V',{\rm Tr}(\F'_k))$ is a BRSC. We claim that $\pure(T_k(S)) = (V',{\rm Tr}(\F'_k))$.

Let $X \in \H \cap P_k(V)$. Then there exists an enumeration $a_1,\ldots,a_k$ of the elements of $X$ such that 
$$\oo{a_1} \subset \oo{a_1a_2} \subset \ldots \subset \oo{a_1\ldots a_k}.$$
Hence $X$ is a transversal of the successive differences for
$$\emptyset \subset \oo{a_1} \subset \oo{a_1a_2} \subset \ldots \subset \oo{a_1\ldots a_{k-1}} \subset V,$$
which is a chain of length $k$ in $\F_k$. Thus $X \in {\rm Tr}(\F'_k)$.

Conversely, assume that $X \in {\rm Tr}(\F'_k)$. We may assume that $X$ is a facet of $(V',{\rm Tr}(\F'_k))$. Then there exists some chain 
\beq
\label{brnm3}
F_0 \subset F_1 \subset \ldots \subset F_s
\eeq
in $\F'_r$ and some enumeration $a_1,\ldots,a_s$ of the elements of $X$ such that $a_i \in F_i \setminus F_{i-1}$ for $i = 1,\ldots,s$.
Since $X$ is a facet, we must have $F_0 = \emptyset$ and $F_s = V$. Suppose that $F_{s-1}\rho = r < k-1$. Since $F_{s-1} \in \F'_k$, then it must occur in some chain of length $k$ in $\F'_k$, hence 
we have some chain 
$$F_{s-1} = F'_0 \subset F'_1 \subset \ldots \subset F'_t \subset F'_{t+1} = V$$ 
in $\F'_k$ for some $t\geq 1$. Since $a_s \in F'_{t+1} \setminus F'_0$, we have $a_s \in F'_{j} \setminus F'_{j-1}$ for some $j \in \{1,\ldots,t+1\}$, hence there exists some $Y \in {\rm Tr}(\F'_k) \cap P_{s+t}(V)$ containing (strictly) $X$, contradicting $X \in \fct(V',{\rm Tr}(\F'_k))$. Thus $F_{s-1}\rho = k-1$. 

Now $a_i \in F_i \setminus F_{i-1}$ for $i = 1,\ldots,s-1$ and so we can apply Lemma \ref{step} $s-1$ times to refine (\ref{brnm3}) to a chain of length $k$ in $L(S)$ of the form
$$F_0 \subset \oo{F_0 \cup a_1} \subseteq \ldots \subseteq F_1 \subset \oo{F_1 \cup a_2} \subseteq \ldots \subseteq F_{s-1} \subset F_s,$$
which admits a transversal of the successive differences containing $X$. Since $X \in \fct(V',{\rm Tr}(\F'_k))$, it follows that $s = k$ and so in view of Lemma \ref{opre} we have $X \in \H \cap P_k(V)$, hence $X$ is a facet of $\pure(T_k(S))$. Therefore $\pure(T_k(S)) = (V',{\rm Tr}(\F'_k))$ as claimed.
\qed

Together with Proposition \ref{bpnm}, this yields:

\bc
\label{kan}
Problems \ref{pc}, \ref{tpc} and \ref{tpct} have positive answers for boolean representable near-matroids. In particular, they hold for:
\bi
\item[(i)] paving BRSC;
\item[(ii)] BRSC of dimension $\leq 2$. 
\ei
\ec

As remarked earlier, Example \ref{pobe} answers negatively Problem \ref{tpc} for dimension 3 and $k = 4$. On the other hand, Problem \ref{tpc} has a positive answer for $k \leq 2$: if $S$ is a BRSC, then $T_k(S)$ is a BRSC by Proposition \ref{tbone} and $\pure(T_k(S))$ is a BRSC by Corollary \ref{kan}(ii).

The next example (discussed in the Appendix) answers negatively Problem \ref{tpc} for $\dim(S) = 3$ and $k = 3$.

\be
\label{cepc}
Let $S = (V,\H)$ with $V = \cup\{ i,i',i'' \mid i \in \Z_3 \}$,
$$Z = \cup\{ i(i+1)(i+1)', i''(i+1)(i+1)' \mid i \in \Z_3 \}$$
and
$$\begin{array}{lll}
\H&=&(P_{\leq 3}(V) \setminus Z) \cup \{ ii''(i+1)p \mid i \in \Z_3,\; p \in V \setminus ii''(i+1)(i+1)' \}\\
&\cup&\{ ii''(i+1)'p \mid i \in \Z_3,\; p \in V \setminus ii''(i+1)(i+1)' \}.
\end{array}$$
Then $S$ is a BRSC but ${\rm pure}(T_3(S))$ is not.
\ee

Another counterexample, also analyzed in the Appendix, is given by the {\em boolean module} $\B^{(4)}$: a simplicial complex of dimension 3 admitting a $4  \times (2^4-1)$ boolean matrix representation where all columns are distinct and nonzero (so we have all possible nonzero columns).

\be
\label{bfour}
The boolean module $\B^{(4)}$ is pure and its truncation to rank 3 is a pure TBRSC which is not a BRSC.
\ee

We turn now our attention to Problem \ref{tpct}.

As remarked earlier, Example \ref{pobe} also answers negatively Problem \ref{tpct} for dimension 3 and $k = 4$.
The next result shows, that, unlike Problem \ref{tpc}, Problem \ref{tpct} admits a positive answer for $k \leq 3$.

\bt
\label{spch}
Let $S$ be a BRSC and let $1 \leq k \leq 3$. Then {\rm pure}$(T_k(S))$ is a TBRSC.
\et

\proof
In the Appendix.
\qed

\section{Topology}

In this section, we generalize to TBRSCs results proved in \cite{MRS} for the topology of BRSCs.

Let $S = (V,\H)$ be a simplicial complex. 
We say that
$S$ is {\em connected} if the graph $T_2(S)$ is connected. The proof for the following result is essentially the proof given for BRSCs in \cite[Lemma 3.1]{MRS}.

\bl
\label{connected}
Let $S = (V,\H)$ be a TBRSC. Then
$S$ is connected unless $\H = P_1(V)$ and $|V| > 1$.
\el

\proof
In the Appendix.
\qed

It is well known that the geometric realization $||S||$ of a simplicial
complex $S$, a subspace of some euclidean space $\mathbb{R}^n$, is unique
up to homeomorphism. For details, see e.g. \cite[Appendix A.5]{RSm}. 

Given a point $v_0 \in ||S||$, the {\em
  fundamental group} $\pi_1(||S||,v_0)$ 
is the group having as
elements the homotopy equivalence classes of closed paths
$$
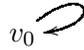
\xymatrix{
v_0 \ar@(ur,r)
}$$
the product being determined by the concatenation of paths. 

If $S$ is connected, then $\pi_1(||S||,v_0) \cong \pi_1(||S||,w_0)$
for all points $v_0,w_0$ in $||S||$, hence we may use the notation
$\pi_1(||S||)$ without ambiguity. We produce now a presentation for
$\pi_1(||S||)$. This combinatorial description is also known as the
{\em edge-path group} of $S$ (for details on the
fundamental group of a simplicial complex, see \cite{Spa}). 

We fix a spanning tree $T$ of $S$ and we define 
$$A = \{ a_{pq} \mid pq \in \H \cap P_2(V)\},$$
$$R_T = \{ a_{qp}a_{pq}\inv \mid pq \in \H \cap P_2(V) \} \cup \{
a_{pq}a_{qr}a_{pr}\inv \mid pqr \in \H \cap P_3(V) \} 
\cup \{ a_{pq} \mid pq \in T \}.$$
We may view $\pi_1(||S||)$ as the group
defined by the group presentation 
\beq
\label{pre}
\langle A \mid R_T \rangle.
\eeq

We compute next the fundamental group of a connected TBRSC. If it has dimension 1, it is a graph and so it follows easily from the presentation (\ref{pre}) that its fundamental group is free of rank $e-v+1$, where $e$ (respectively $v$ denotes the number of edges (respectively vertices). Note that $v-1$ is the number of edges of the spanning tree $T$. Therefore we concentrate our attention in the case of dimension $\geq 2$. These TBRSCs are connected by Lemma \ref{connected}.
 
Given a BRSC $S = (V,\H)$, the {\em graph of flats} $\Gamma(L(S))$ has
vertex set $V$ and edges $p \edge q$ whenever $p \neq q$ and $\oo{pq}
\subset V$. 

Let $C$ be a connected component of $\Gamma(L(S))$. If $\H \cap
P_2(C) \neq \emptyset$, we shall say that $C$ is 
$\H$-{\em nontrivial}. Otherwise, we say that $C$ is $\H$-{\em trivial}.
The {\em size} of $C$ is its number of vertices.

The next result shows that, given a TBRSC $S = (V,\H)$ of dimension $\geq 2$, the graph of flats $\Gamma(L(S^{\varepsilon}))$ and the size of its
$\H^{\varepsilon}$-trivial components determine completely
the fundamental group of $S$. Note that $L(S^{\varepsilon}) = \varepsilon(S)$ by Theorem \ref{eqtr}, hence, for all distinct $p,q \in V$, $p \edge q$ is an edge of $L(S^{\varepsilon}) $ if and only if there exists some $Z \in \varepsilon(S)$ such that $pq \subseteq Z \subset V$.

\bt
\label{fung}
Let $S$ be a TBRSC of dimension $\geq 2$. Assume that $\Gamma(L(S^{\varepsilon}))$ has $s$
$\H^{\varepsilon}$-nontrivial connected components and $r$ $\H^{\varepsilon}$-trivial connected
components of sizes $f_1,\ldots,f_r$. 
Then $\pi_1(||\S||)$ is a free group of rank 
$$\binom{s+f_1+\ldots +f_r-1}{2} - \sum_{i=1}^r \binom{f_i}{2},$$
or equivalently, 
$$\binom{s-1}{2} +(s-1)(f_1+\ldots +f_r) + \sum_{1 \leq i < j \leq r} f_if_j.$$
\et

\proof
This result was proved in \cite[Theorem 3.3]{MRS} for BRSCs (with $S^{\varepsilon}$ replaced by $S$). Therefore it suffices to note that $S$ and $\S^{\varepsilon}$ have the same fundamental group. Indeed, $\pi(||\S||) = \pi(||T_3(S)||)$ and $\pi(||\S^{\varepsilon}||) = \pi(||T_3(S^{\varepsilon})||)$. Since $S$ has dimension $\geq 2$, it follows from Theorem \ref{eqtr} that $T_3(S) = T_3(S^{\varepsilon})$. It follows that $\pi(||S||) = \pi(||S^{\varepsilon}||)$ as required.
\qed

\bc
\label{simfun}
Let $S$ be a simple TBRSC of dimension $\geq 2$. 
Then $\pi_1(||S||)$ is a free group of rank $\binom{t-1}{2}$, where
$t$ denotes the number of connected components of $\Gamma(L(S^{\varepsilon}))$. 
\ec

\proof
If $S = (V,\H)$ is simple, then each $\H^{\varepsilon}$-trivial connected component of
$\Gamma(L(S^{\varepsilon}))$ has precisely one vertex. Hence, by Theorem
\ref{fung}, $\pi_1(||S||)$ is a free group of rank
$\binom{t-1}{2}$.
\qed

In \cite[Example 3.5]{MRS}, it is shown that free groups of rank $\binom{n}{2}$ $(n \geq 2)$
occur effectively 
as fundamental groups of simple BRSCs of dimension 2.

Let $S = (V, \J)$ be a simplicial complex. We recall now the definitions
of the (reduced) homology groups of $S$ (see e.g. \cite{Hat}).

If $S$ has $c$ connected components, it is well known that the {\em
  0th homology group} $H_0(S)$ is isomorphic to 
the free abelian group of rank $c$. For dimension $k \geq 1$, we
proceed as follows.

Fix a total ordering of $V$.
Let $C_k(S)$ denote the free abelian group on $\J \cap P_{k+1}(V)$,
that is, all the formal sums of the 
form $\sum_{i \in I} 
n_iX_i$ with $n_i \in \Z$ and $X_i \in \J \cap P_{k+1}(V)$
(distinct). Given $X \in \J \cap P_{k+1}(V)$, write $X = x_0x_1\ldots x_k$
with $x_0 < \ldots < x_k$. We define 
$$X\partial_k = \sum_{i= 0}^k (-1)^i (X\setminus\{ x_i\}) \in
C_{k-1}(S)$$
and extend this by linearity to a homomorphism $\partial_k:C_{k}(S)
\to C_{k-1}(S)$ (the $k$th {\em boundary map} of $S$). Then the {\em
  $k$th homology group} of $S$ is defined as the quotient
$$H_k(S) = \ker\partial_k / {\rm Im}\partial_{k+1}.$$

The {\em $0$th reduced homology group} of $S$, denoted by
$\tilde{H}_0(S)$, is isomorphic to the free abelian group of rank $c-1$, where
$c$ denotes the number of connected components of $S$.
For $k \geq 1$, the {\em $k$th reduced homology group} of $S$,
denoted by $\tilde{H}_k(S)$ coincides with the $k$th homology group.

A {\em wedge} of spheres $S_1,\ldots, S_m$ (of possibly different
dimensions) is a 
topological space obtained by identifying $m$ points $s_i \in S_i$ for
$i = 1,\ldots,m$.  

We say that two topological spaces $X$ and $Y$ have the same {\em homotopy
  type}
if there exist continuous mappings $\alpha:X \to Y$ and $\beta:Y \to
X$ such that:
\bi
\item
there exists a homotopy between $\alpha\beta$ and $1_X$;
\item
there exists a homotopy between $\beta\alpha$ and $1_Y$.
\ei
An important theorem of Bj\"orner and Wachs \cite{BW} states that shellable simplicial complexes (a class including matroids as a particular case)  have the homotopy type of a wedge of spheres.

Theorem \ref{fung} also yields the following important consequence, where the proof is essentially the proof given for BRSCs in \cite[Theorem 3.6]{MRS}.

\bt
\label{wed}
Let $S$ be a TBRSC of dimension $2$. Then:
\bi
\item[(i)] the homology groups of $S$ are free abelian;
\item[(ii)] $S$ has the homotopy type of a wedge of 1-spheres and
  2-spheres.
\ei
\et

\proof
In the Appendix.
\qed

\section*{Acknowledgments}

The first author acknowledges support from the Binational Science Foundation (BSF) of the United States and Israel, grant number 2012080. The second author acknowledges support from the Simons Foundation (Simons Travel Grant Number 313548).
The third author was partially supported by CMUP (UID/MAT/00144/2019), which is funded by FCT (Portugal) with national (MCTES) and European structural funds through the programs FEDER, under the partnership agreement PT2020.

% ----------------------------------------------------------------
%\bibliographystyle{amsplain}

\section*{Appendix}

We collect in this Appendix several proofs omitted from the main text, and the discussion of several examples.

\medskip
\noindent
{\em Proof of Proposition \ref{tru}}. Let $X \in
L(T_k(S))$. If $X$ contains a facet of $T_k(S)$, then $X = V$, so we assume that $X$ contains no facet of $T_k(S)$. Now let $I \in \H \cap 2^X$ and $p \in V
\setminus X$. Since $I \notin \fct(T_k(S))$, we have $|I| < k$ and so $I
\in T_k(\H)$. Now $X \in L(T_k(S))$ yields $I \cup \{ p \} \in T_k(\H) \subseteq \H$. Therefore $X \in
L(S)$ and the direct inclusion holds.

Conversely, assume that $X \in
L(S)$ contains no facet of $T_k(S)$. Let $I \in \H \cap P_{\leq k}(X)$ and $p \in V \setminus X$. Since
$X \in L(S)$, we get $I \cup \{ p \} \in
\H$. But $I$ is not a facet of 
$T_k(S)$, hence $|I| < k$ and so $I \cup \{ p \} \in T_k(\H)$. Thus $X \in
L(T_k(S))$ as required.
\qed

\noindent
{\em Proof of Lemma \ref{inc}}. 
(i)
We prove that
\beq
\label{inc1} 
T_k(\H^{\varepsilon}) \subseteq \H
\eeq
holds for $k = 0,\ldots, d+1$ by induction on $k$.

The case $k = 1$ being trivial, assume that $k \in \{ 2, \ldots, d+1\}$ and
(\ref{inc1}) holds for $k-1$. Let $X \in T_k(\H^{\varepsilon})$. We may assume that $|X| = k$. 
Then there exists an enumeration 
$x_1,\ldots, x_k$ of $X$ and $Z_0, \ldots, Z_k \in \varepsilon(S)$ such that
$Z_0 \supset Z_1 \supset \ldots \supset Z_k$ and $x_i \in Z_{i-1}
\setminus Z_i$ for $i = 1, \ldots, k$. Let $X' = \{ x_2, \ldots, x_k
\}$. Since $X' \in T_{k-1}(\H^{\varepsilon})$, it follows from the
induction hypothesis that $X' \in \H$. Now $|X'| \leq d$, $X' \subseteq
Z_1$ and $x_1 \in V \setminus Z_1$, hence it follows from $Z_1 \in
\varepsilon(S)$ that $X = X' \cup \{ x_1 \} \in \H$. Thus (\ref{inc1}) holds
for $k = 1,\ldots, d+1$.

(ii) Let $X \in \varepsilon(S)$. Let $I \in \H^{\varepsilon} \cap 2^X$ and $p \in V
\setminus X$. Since $I \in \H^{\varepsilon}$, there exists an enumeration
$x_1,\ldots, x_k$ of $I$ and $Z_0, \ldots, Z_k \in \varepsilon(S)$ such that
$Z_0 \supset Z_1 \supset \ldots \supset Z_k$ and $x_i \in Z_{i-1}
\setminus Z_i$ for $i = 1, \ldots, k$. Now by Lemma \ref{propt}(i)
$$Z_0\cap X \supset Z_1 \cap X  \supset \ldots \supset Z_k\cap X$$
is also a chain in $\varepsilon(S)$ satisfying $x_i \in (Z_{i-1} \cap X)
\setminus (Z_i \cap X)$ for $i = 1, \ldots, k$. Since $V \supset Z_0
\cap X$ is also a chain in $\varepsilon(S)$ and $p \in V \setminus (Z_0 \cap
X)$, we get $I \cup \{ p \} \in \H^{\varepsilon}$ and so $X \in
L(S^{\varepsilon})$.

(iii) Let $X \in \H^{\varepsilon}$. Then $X$ 
is a transversal of the partition of successive differences for some
chain of $\varepsilon(S)$, and so is any subset of $X$. Thus $S^{\varepsilon}$ is a simplicial complex. By (ii), a chain in $\varepsilon(S)$ is also a chain in $L(S^{\varepsilon})$. Therefore $S^{\varepsilon}$ is boolean representable.
\qed

\noindent
{\em Proof of Theorem \ref{eqtr}}. 
(i) $\Rw$ (ii). Write $S' = (V,\H')$. We start by showing that 
\beq
\label{eqtr1}
L(S') \subseteq \varepsilon(S).
\eeq
Let $F \in L(S')$. Suppose that $X \in \H \cap P_{\leq d}(F)$ and $p \in V
\setminus F$. Since $\H \subseteq \H'$, it follows from $F \in
L(S')$ that $X \cup \{ p \} \in \H'$. But now $|X| \leq d$ implies $X
\cup \{ p \} \in T_{d+1}(\H') = \H$ and so $F \in \varepsilon(S)$. Therefore
(\ref{eqtr1}) holds.

Now let $X \in \H$. Since $\H \subseteq \H'$, there exists an enumeration
$x_1,\ldots, x_k$ of $X$ and $F_0, \ldots, F_k \in L(S')$ such that
$F_0 \supset F_1 \supset \ldots \supset F_k$ and $x_i \in F_{i-1}
\setminus F_i$ for $i = 1, \ldots, k$. By (\ref{eqtr1}), we have $F_0,
\ldots, F_k \in \varepsilon(S)$ and so $X \in \H^{\varepsilon}$. Since $\dim(S) = d$,
then $X \in T_{d+1}(\H^{\varepsilon})$ and so $\H \subseteq T_{d+1}(\H^{\varepsilon})$. Therefore $\H =
T_{d+1}(\H^{\varepsilon})$ by Lemma \ref{inc}(i), and so $S = T_{d+1}(S^{\varepsilon})$.

(ii) $\Rw$ (i). This follows from Lemma \ref{inc}(iii). 

It remains to be proved that $L(S^{\varepsilon}) = \varepsilon(S)$. 

Let $X \in L(S^{\varepsilon})$. Let $I \in \H \cap P_{\leq d}(X)$ and $p \in V
\setminus X$. Then $I \in \H^{\varepsilon}$ by (ii) and so $X \in L(S^{\varepsilon})$
yields $I \cup \{ p \} \in \H^{\varepsilon}$. Since $|I| \leq d$, we get $I \cup \{
p \} \in T_{d+1}(\H^{\varepsilon}) = \H$ and so $X \in \varepsilon(S)$. The opposite inclusion
follows from Lemma \ref{inc}(ii).
\qed

\noindent
{\em Analysis of Example \ref{truno}}. 
Indeed, it is easy to check that 
$$\emptyset \subset 1 \subset 12 \subset 1235 \subset V, \quad 
\emptyset \subset 3 \subset 1235 \subset V, \quad 
\emptyset \subset 4 \subset V$$
are all 
chains in $\varepsilon(S)$. Now every $X \in \H$ is a partial transversal of either chain (if $ X = 46$, we use the third 
chain, if $X \supseteq 35$ we use the second chain, 
in the remaining cases we use the first). Hence $\H \subseteq \H^{\varepsilon}$ and so $S = T_3(S^{\varepsilon})$ by Lemma \ref{inc}(i). Therefore $S$ is a TBRSC by Theorem \ref{eqtr}.

Consider now $134 \in \H$. 
\bi
\item
Since $135 \notin \H$, we get $5 \in \oo{13}$. Since $235 \notin \H$, we get $2 \in \oo{35} \subseteq \oo{13}$. Since $123 \in \fct(S)$, we get $\oo{13} = V$.
\item
Since $146 \notin \H$, we get $6 \in \oo{14}$. Since $246 \notin \H$, we get $2 \in \oo{46} \subseteq \oo{14}$. Since $124 \in \fct(S)$, we get $\oo{14} = V$.
\item
Since $346 \notin \H$, we get $6 \in \oo{34}$. Since $246 \notin \H$, we get $2 \in \oo{46} \subseteq \oo{34}$. Since $234 \in \fct(S)$, we get $\oo{34} = V$.
\ei
It follows that $S$ is not a BRSC.
\qed
%
%OR:
%
%Indeed, it is easy to check that $\flatx\H$ is the lattice
%$$\xymatrix{
%&&V \ar@{-}[dl] &&&\\
%&12 &&&&\\
%1 \ar@{-}[ur] & 2 \ar@{-}[u] & 3 \ar@{-}[uu]
%& 4 \ar@{-}[uul] & 5 \ar@{-}[uull] & 6
%\ar@{-}[uulll] \\
%&&\emptyset \ar@{-}[ull] \ar@{-}[ul] \ar@{-}[u] \ar@{-}[ur] \ar@{-}[urr]
%\ar@{-}[urrr] &&&
%}$$
%and so $\H \notin \bpav(2)$ since 135 is not a transversal of any chain in $\flatx\H$.
%
%On the other hand, $T(H)$ is the lattice
%$$\xymatrix{
%&&&V \ar@{-}[dl] &&\\
%&&1234 \ar@{-}[dl] &&&\\
%&12 &&&&\\
%1 \ar@{-}[ur] & 2 \ar@{-}[u] & 3 \ar@{-}[uu]
%& 4 \ar@{-}[uul] & 5 \ar@{-}[uuul] & 6
%\ar@{-}[uuull] \\
%&&\emptyset \ar@{-}[ull] \ar@{-}[ul] \ar@{-}[u] \ar@{-}[ur] \ar@{-}[urr]
%\ar@{-}[urrr] &&&
%}$$
%and it is easy to check that $\H = (\J(T(H)))_{3}$. Therefore $\H \in \tbpav(2)$.

\noindent
{\em Analysis of Example \ref{extruA}}. Let $a,b,c$ be an enumeration of $123$. Let $X \in \varepsilon(S)$ contain $ab$. Since $ab4 \notin \H$, we have $4 \in X$. Since $ac4 \notin \H$, we get $c \in X$. Hence $X = V$ and so 123 cannot be a transversal of the successive differences for a chain in $\varepsilon(S)$. Therefore $S$ is not a TBRSC.
\qed 

\noindent
{\em Proof of Proposition \ref{eqp}}. 
Let $F \in L(S^{\varepsilon})$. Take $X \in \H \cap P_{\leq d}(F)$ and $p \in V \setminus F$. Since $P_{\leq d}(V) \subseteq \H$, we have $P_{\leq d-1}(V) \subseteq L(S)$, hence
$P_{\leq d-1}(V) \subseteq \varepsilon(S) \subseteq L(S^{\varepsilon})$ by Lemmas \ref{propt}(ii) and \ref{inc}(ii). Thus $X \in \H^{\varepsilon}$ and since $F \in L(S^{\varepsilon})$ we get $X \cup \{ p \} \in \H^{\varepsilon}$. Thus $X \cup \{ p \} \in T_{d+1}(\H^{\varepsilon}) \subseteq \H$ by Lemma \ref{inc}(i) and so $F \in \varepsilon(S)$.
\qed

\noindent
{\em Analysis of Example \ref{ttnot}}. 
Straightforward computation shows that $\varepsilon(S) = \{ \emptyset, 35, V\}$. It follows easily that $124 \in L(S^{\varepsilon}) \setminus \varepsilon(S)$.
\qed

\noindent
{\em Proof of Lemma \ref{buildm}}. 
It is immediate that $(V,\H)$ is a simplicial complex.
Let $I,J \in \H$ with $|I| = |J|+1$. We may assume that $|J \cap I| \leq 1$.

Assume first that $J \cap I = \{ a\}$. Write $J = ab$ and $I = ac_1c_2$. Suppose that $abc_s \notin \H$ for $s = 1,2$. Then $abc_s$ contains some $F_s \in \F$ for $s = 1,2$. Since $I,J \in \H$, we must have $F_s = bc_s$ for $s = 1,2$. But then $F_1 \cap F_2 = \{ b\}$, a contradiction. Thus $J \cup \{c_s\} \in \H$ for some $s$.

Assume now that $J \cap I = \emptyset$. Write $J = ab$ and $I = c_1c_2c_3$. Suppose that $abc_s \notin \H$ for $s = 1,2,3$. Then $abc_s$ contains some $F_s \in \F$ for $s = 1,2,3$. Since $J \in \H$, we must have $F_s \in \{ ac_s, bc_s\}$ for $s = 1,2,3$. But then there exist $i,j \in \{ 1,2,3\}$ such that $|F_i \cap F_j| = 1$, a contradiction. Thus $J \cup \{c_s\} \in \H$ for some $s$.
\qed

\noindent
{\em Analysis of Example \ref{nonun}}. 
Indeed, $S_1$ is a uniform matroid and $S_2$ is a matroid by Lemma \ref{buildm}. We may write $S_1 \vee S_2 = (V,\H)$ with 
$$\H = P_{\leq 2}(V) \cup \{ X \in P_3(V) \mid, 12,34 \not\subseteq X \}.$$
We have $1235 \in \H$. Let $Z \in \varepsilon(\H)$.

If $13 \subseteq Z$, then $123 \notin \H$ yields $2 \in Z$, and $125 \notin \H$ yields $5 \in Z$.

If $15 \subseteq Z$, then $125 \notin \H$ yields $2 \in Z$, and $123 \notin \H$ yields $3 \in Z$. Out of symmetry, $35 \subseteq Z$ implies $1 \in Z$.

It follows that $135 \notin \H^{\varepsilon}$ and so $S_1 \vee S_2$ is not a TBRSC by Theorem \ref{eqtr}.
\qed

\noindent
{\em Analysis of Example \ref{ncu}}. 
Indeed, it is easy to check that
$$L(V,\H) \supseteq P_{\leq 1}(V) \cup \{ 12, V\}, \quad
L(V,\H') \supseteq P_{\leq 1}(V) \cup \{ 1235, V\},$$
and it follows easily that 
$(V,\H),(V,\H') \in {\rm BPav}(d)$. We have seen in Example \ref{truno} that $(V,\H\cup \H') \notin {\rm BPav}(d)$.
\qed

\noindent
{\em Analysis of Example \ref{scznp}}.
First, we show that $S$ is not a TBRSC. Suppose it is. Then $1246 \in \H \subseteq \H^{\varepsilon}$ by Theorem \ref{eqtr}. Then there exists some $Z \in \varepsilon(S)$ such that $|1246 \cap Z| = 3$.  Out of symmetry, we may assume that $24 \subset Z$.
\bi
\item
If $1246 \cap Z = 124$, then $124  \in \H \cap 2^Z$ and $1234 \notin \H$ yield $3 \in Z$. Now 
$123  \in \H \cap 2^Z$ and $1236 \notin \H$ yield $6 \in Z$, a contradiction.
\item
If $1246 \cap Z = 246$, then $246  \in \H \cap 2^Z$ and $2346 \notin \H$ yield $3 \in Z$. Now 
$234  \in \H \cap 2^Z$ and $1234 \notin \H$ yield $1 \in Z$, also a contradiction.
\ei
Therefore $S$ is not a TBRSC.

Now let $\mathcal{J} \subseteq 2^V$ be the set of partial transversals of the partial differences for the chain 
$$\emptyset \subset 1 \subset 123 \subset 12345 \subset V.$$
Then $(V,\mathcal{J})$ is a BRSC and $\mathcal{J} \subseteq \H$. On the other hand, $P_{\leq 3}(V) \subseteq \H$. Since $\mathcal{J} \cup P_{\leq 3}(V) = \H$ and $S$ is not a TBRSC, it follows that $S$ admits no largest truncated boolean representable subcomplex.
\qed

\noindent
{\em Proof of Proposition \ref{six}}. 
We fix $V = \{ 1,\ldots,6\}$ as the set of points and we consider $S = (V,\H) \in \tbpav(2) \setminus \bpav(2)$. Then there exists some BRSC $S' = (V,\H')$ such that $S = T_3(S')$. Given $X \subseteq V$, let $\oo{X}$ (respectively $\wh{X}$) denote the closure of $X$ in $L(S')$ (respectively $L(S)$). 

Since $S \notin \bpav(2)$ and $P_{\leq 1}(V) \subseteq L(S)$, there exists some $X \in P_3(V) \cap \H$ such that 
\begin{equation}
\label{six2}
X \subseteq \wh{X \setminus \{x\}} \mbox{ for every }x \in X.
\end{equation} 
Without loss of generality, we may assume that $X = 345$. On the other hand, since $S' \in \bpav(2)$ and $345 \in \H'$, there exists some $x \in 345$ such that $x \notin \oo{345 \setminus \{x\}}$. We may assume that $x = 5$. We claim that
\begin{equation}
\label{six1}
|\oo{34}| = 4.
\end{equation}
Indeed, we know already that $5 \notin \oo{34}$. Suppose that $\oo{34} = 34$. Then $34y \in \H'$ (and therefore $34y \in \H$) for every $y \in 1256$, yielding $\wh{34} = 34$, contradicting (\ref{six2}). Without loss of generality, we may assume that $34y \notin \H'$ for some $y \in 126$, say $y = 1$. Hence $134 \subseteq \oo{34}$. Suppose that $\oo{34} = 134$. Since $134 \notin \H$, this implies $\wh{34} = \oo{34} = 134$, contradicting (\ref{six2}). Thus $|\oo{34}| \geq 4$. Since $5 \notin \oo{34}$, we may assume without loss of generality that $1234 \subseteq \oo{34}$. 

Suppose that $1234 \subset \oo{34}$. Since $5 \notin \oo{34}$, we get $\oo{34} = 12346$. It follows that $45z \in \H'$ for every $z \in 1236$, hence $\wh{45} = 45$, contradicting (\ref{six2}). Therefore $\oo{34} = 1234$ and so (\ref{six1}) holds. 

It follows that $ab5, ab6 \in \H$ for all $a,b \in 1234$ distinct. Since $\oo{134} = 1234$, it follows that $\{ 123, 124, 234\} \not\subseteq \H$. Together with $134 \notin \H$, this implies that the restriction 
$$S'' =  S'|_{1234} = S|_{1234} = (1234,\H'')$$
misses at least two triangles. 

On the one hand, $134 \notin \H$ and $\{ 123, 124, 234\} \not\subseteq \H$ yield $1234 \subseteq \wh{34}$. 
On the other hand, it follows from (\ref{six2}) that $5 \in \wh{34}$, hence $12345 \subseteq \wh{34}$. Since $ab5, ab6 \in \H'$ for all $a,b \in 1234$ distinct, then $1234 \notin L(S)$ (if $1234 \in L(S)$, then $\wh{34} \subseteq 1234$)
implies that $1234 \setminus \{c\} \in \H$ for some $c \in 1234$. Therefore $S''$ has exactly one or two triangles. Since $S''$ is a restriction of the BRSC $S'$, it follows that $S''$ is a BRSC. 
On the other hand, it follows from \cite[Example 5.2.11]{RSm} that a paving BRSC with 4 points cannot have exactly one triangle, hence $S''$ has exactly two triangles, whose intersection has two points, say $de$. 

Together with $1234 \in L(S')'$, this implies that $de \in L(S')$. Since $134 \notin \H$, we have $de \in \{ 12,23,24\}$. Since we have not distinguished 3 from 4 so far, we may assume that $de \in \{ 12,23\}$. 

In any case, having $1234 \in L(S')$ determines that $ab5, ab6 \in \H$ for all $a,b \in 1234$ distinct (12 elements), and $de \in L(S')$ determines which two elements among the four elements of $P_3(1234)$ belong to $\H$. Thus we only need to discuss what happens with $156,256,356,456$. If $356 \in \H$, then $35 \in L(S')$ (in view of $1234 \in L(S')$), implying $\oo{35} = 35$ (and consequently $\wh{35} = 35$), contradicting (\ref{six2}). Therefore $356 \notin \H$. Similarly, $456 \notin \H$. It follows that we reduced the discussion to determine whether or not  $156, 256 \in \H$, for each choice of $de \in \{ 12,23\}$.

If we omit both $156, 256$ from $\H$, we get the two cases
\bi
\item[(1)] $\H = \B_2(1234) \cup \B_2(12)$, 
\item[(1')] $\H = \B_2(1234) \cup \B_2(23)$,
\ei
which are clearly isomorphic.

Now adding $156$ (respectively $256$) is the only consequence of adding $15$ (respectively $25$) as a line, and these additions do not interfere with each other. We are then bound to consider the cases: 
\bi
\item[(2)] $\H = \B_2(1234) \cup \B_2(12) \cup \B_2(15)$;
\item[(2')] $\H = \B_2(1234) \cup \B_2(12) \cup \B_2(25)$;
\item[(3)] $\H = \B_2(1234) \cup \B_2(12) \cup \B_2(15) \cup \B_2(25)$;
\item[(4)] $\H = \B_2(1234) \cup \B_2(23) \cup \B_2(15)$;
\item[(2'')] $\H = \B_2(1234) \cup \B_2(23) \cup \B_2(25)$;
\item[(5)] $\H = \B_2(1234) \cup \B_2(23) \cup \B_2(15) \cup \B_2(25)$. 
\ei

The cases (2), (2') and (2'') are clearly isomorphic. Applying the permutations $(13)$ and $(132)$ to $12345$ in cases (4) and (5), respectively, we have reduced our discussion to the cases
\bi
\item[(1)] $\H = \B_2(1234) \cup \B_2(12)$;
\item[(2)] $\H = \B_2(1234) \cup \B_2(12) \cup \B_2(15)$;
\item[(3)] $\H = \B_2(1234) \cup \B_2(12) \cup \B_2(15) \cup \B_2(25)$;
\item[(4)] $\H = \B_2(1234) \cup \B_2(12) \cup \B_2(35)$;
\item[(5)] $\H = \B_2(1234) \cup \B_2(12) \cup \B_2(15) \cup \B_2(35)$. 
\ei

We list below the triangles missing in each of the cases:
\bi
\item[(1)] 134, 234, 156, 256, 356, 456;
\item[(2)] 134, 234, 256, 356, 456;
\item[(3)] 134, 234, 356, 456;
\item[(4)] 134, 234, 156, 256, 456;
\item[(5)] 134, 234, 256, 456.
\ei

Out of cardinality arguments, we only have to distinguish (2) from (4) and (3) from (5).
Now 1 appears only once among the missing triangles in (2), and all points appear more often in (4); 1 and 2 appear only once among the missing triangles in (3), but only 1 has a single occurrence in (5). Therefore these complexes (1) -- (5) are nonisomorphic.

By construction, any one of these 5 complexes is in $\tbpav(2)$. We confirm now that neither of them is a BRSC. For the first three cases, we take $345 \in \H$.

\bi
\item[(1)] $134 \notin \H$, hence $1 \in \wh{34}$; $234 \notin \H$, hence $2 \in \wh{34}$; $\wh{34}$ contains the facet $123$, hence $\wh{34} = V$. $356 \notin \H$, hence $6 \in \wh{35}$; $456 \notin \H$, hence $4 \in \wh{35}$. Similarly, $3 \in \wh{45}$.
\item[(2)] Same argument as in (1).
\item[(3)] Same argument as in (1).
\ei

For the remaining two cases, we take $245 \in \H$.
\bi
\item[(4)] $234 \notin \H$, hence $3 \in \wh{24}$; $134 \notin \H$, hence $1 \in \wh{24}$; $\wh{24}$ contains the facet $123$, hence $\wh{24} = V$. $256 \notin \H$, hence $6 \in \wh{25}$; $456 \notin \H$, hence $4 \in \wh{25}$. Similarly, $2 \in \wh{45}$.
\item[(5)] Same argument as in (4).
\ei
\qed

\noindent
{\em Proof of Lemma \ref{pret}}.
Let $S = (V,\H) \in \tbr$ and let $\emptyset \neq W \subseteq V$. Since $S \in \tbr$, there exist a BRSC $S' = (V,\H')$ and $m \geq 1$ such that $S = T_m(S')$. We claim that
\beq
\label{pret1}
S|_W = T_m(S'|_W).
\eeq
This is equivalent to the equality 
\beq
\label{pret2}
\H \cap 2^W = (\H' \cap 2^W) \cap P_{\leq m}(W).
\eeq
Now $S = T_m(S')$ yields $\H = \H' \cap P_{\leq m}(V)$ and so
$$\H \cap 2^W = (\H' \cap P_{\leq m}(V)) \cap 2^W = (\H' \cap 2^W) \cap P_{\leq m}(W).$$
Hence (\ref{pret2}) and consequently (\ref{pret1}) do hold.

Since BRSCs are closed under restriction, then $S'|_W$ is a BRSC and it follows from (\ref{pret1}) that $S|_W \in \tbr$. Thus $\tbr$ is closed under restriction. Since it is also closed under isomorphism, then $\tbr$ is a prevariety of simplicial complexes.

On the other hand, the class of all finite paving simplicial complexes is a prevariety in view of \cite[Proposition 8.3.1(ii)]{RSm}. Since the intersection of two prevarieties is obviously a prevariety, it follows that $\tbp$ is a prevariety itself.
\qed

\noindent
{\em Analysis of Example \ref{pobe}}.
Since $R$ is a Moore family, $S$ is a BRSC. The maximal chains in $R$ are
\beq
\label{mic1}
\emptyset \subset 1 \subset 12 \subset 1234 \subset V, \quad  
\emptyset \subset 5 \subset 56 \subset 356 \subset V,\eeq
\beq
\label{mic2}
\emptyset \subset 3 \subset 1234 \subset V, \quad  
\emptyset \subset 3 \subset 356 \subset V.
\eeq
Hence $\dim(S) = 3$. Since ${\rm Tr}(R)$ is the set of partial transversals of the successive differences for some of these chains, it follows easily that
$$\begin{array}{lll}
{\rm Tr}(R)&=&(P_{\leq 2}(V) \setminus \{ 134, 157, 167, 234, 257, 267, 457, 467\})\\
&\cup&\{ 123a \mid a \in 567 \} \cup \{ 124a \mid a \in 567 \} \cup \{ 356b \mid b \in 1247 \}.
\end{array}.$$
Write ${\rm pure}(S) = (V,\H')$. It is routine to check that 
$$\H' = {\rm Tr}(R) \setminus \{ 347 \}.$$
Indeed, it is easy to see that each $X \in P_2(V)$ is a partial transversal of the successive differences for some chain of type (\ref{mic1}), and to check which transversals of the successive differences for some chain of type (\ref{mic2}) cannot be obtained through chains of type (\ref{mic1}).

Now we have $1235 \in \H'$. Let $Z \in \varepsilon(\H')$ be such that $|Z \cap 1235| \geq 3$. We show that $1235 \subseteq Z$.
\bi
\item
Suppose that $123 \subseteq Z$. Since $123 \in \H'$ and $1234 \notin \H'$, we have $4 \in Z$. Since $34 \in \H'$ and $347 \notin \H'$, we have $7 \in Z$. Since $47 \in \H'$ and $457 \notin \H'$, we get $5 \in Z$.
\item
Suppose that $125 \subseteq Z$. Since $125 \in \H'$ and $1257 \notin \H'$, we have $7 \in Z$. Since $127 \in \H'$ and $1267 \notin \H'$, we have $6 \in Z$. Since $567 \in \H'$ and $4567 \notin \H'$, we have $4 \in Z$. Since $24 \in \H'$ and $234 \notin \H'$, we get $3 \in Z$.
\item
Suppose that $135 \subseteq Z$. Since $135 \in \H'$ and $1345 \notin \H'$, we have $4 \in Z$. Since $345 \in \H'$ and $2345 \notin \H'$, we get $2 \in Z$.
\item
Suppose that $235 \subseteq Z$. Since $235 \in \H'$ and $2345 \notin \H'$, we have $4 \in Z$. Since $345 \in \H'$ and $1345 \notin \H'$, we get $1 \in Z$.
\ei
Thus there exists no $Z \in \varepsilon(\H')$ such that $|Z \cap 1235| = 3$. By Theorem \ref{eqtr}, ${\rm pure}(S)$ is not a TBRSC.
\qed

\noindent
{\em Proof of Lemma \ref{wool}}.
Let $I' \in \H$ be maximal with respect to $I \subseteq I' \subseteq \oo{J}$. If $\oo{I'} \subset \oo{J}$, we can take $p \in \oo{J} \setminus \oo{I'}$ and get $I' \cup \{ p\} \in \H \cap 2^{\oo{J}}$, contradicting the maximality of $I'$. Thus $\oo{I'} = \oo{J}$ and we are done.
\qed

\noindent
{\em Proof of Lemma \ref{opre}}.
Suppose that $F\rho \geq F'\rho$. Then there exist $I,J \in \H$ such that $F = \oo{I}$, $F' = \oo{J}$ and $|I| \geq |J|$. Hence $I \subseteq \oo{J}$ and so by Lemma \ref{wool} there exists some $I' \in \H$ such that $I \subseteq I'$ and $\oo{I'} = \oo{J}$. But we have then $|I'| > |I| \geq |J|$, a contradiction since $S$ is a near-matroid. Therefore $F\rho < F'\rho$.
\qed

\noindent
{\em Proof of Lemma \ref{step}}.
Write $F = \oo{I}$ with $I \in \H$. Since $a_1 \in F'\setminus F$, we have $I \cup a_1 \in \H$. Thus 
$$F \subset \oo{I\cup a_1} = \oo{F\cup a_1} \subseteq F'.$$
Moreover, $$\oo{F\cup a_1}\rho = |I \cup a_1| = |I|+1 = F\rho+1.$$
If $\oo{F\cup a_1} = F'$, we can now iterate this argument to produce a chain
$$F \subset \oo{F\cup a_1} \subset \oo{F\cup a_1a_2} \subset \ldots \subset \oo{F\cup a_1\ldots a_s} = F'$$
for some $a_2,\ldots a_s \in V$ such that 
$\oo{F\cup a_1\ldots a_j}\rho = \oo{F\cup a_1\ldots a_{j-1}}\rho+1$ for $j = 1,\ldots,s$. Thus $s = F'\rho -F\rho = k$ and we are done.
\qed  

\noindent
{\em Analysis of Example \ref{cepc}}.
It is easy to check that $S$ is indeed a simplicial complex. Clearly, $P_{\leq 1}(V) \subset L(S)$. If $X \in P_2(V)$ is not contained in any element of $Z$, then $\oo{X} = X$. Hence, if $abc \in \H$ and $ab$ is not contained in any element of $Z$, then $abc$ is a transversal of the successive differences for the chain
$$\emptyset \subset a \subset ab \subset V$$
in $L(S)$. On the other hand, it is easy to check that the unique $X \in P_3(V) \cap \H$ having all 2-subsets contained in elements of $Z$ is $123$ (see the picture below, where the yellow triangles are the elements of $Z$):

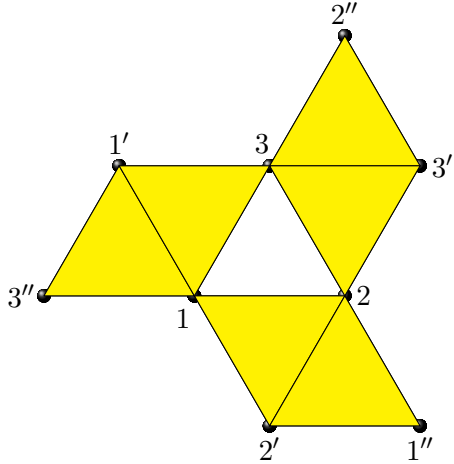
\begin{center}
\begin{tikzpicture}[domain=-4:4]
\draw (0,0) -- (4,0);
\draw (1,1.73) -- (5,1.73);
\draw (3,-1.73) -- (5,-1.73);
\draw (0,0) -- (1,1.73);
\draw (2,0) -- (4,3.46);
\draw (3,-1.73) -- (5,1.73);
\draw (1,1.73) -- (3,-1.73);
\draw (3,1.73) -- (5,-1.73);
\draw (5,1.73) -- (4,3.46);
\node at (1,2.03) {$1'$};
\node at (4,3.76) {$2''$};
\node at (2.9,2.03) {$3$};
\node at (5.3,1.73) {$3'$};
\node at (1.85,-0.3) {$1$};
\node at (-0.3,0) {$3''$};
\node at (4.25,0) {$2$};
\node at (3,-2.03) {$2'$};
\node at (5,-2.03) {$1''$};
\shade[shading=ball, ball color=black] (1,1.73) circle (.09);
\shade[shading=ball, ball color=black] (3,1.73) circle (.09);
\shade[shading=ball, ball color=black] (5,1.73) circle (.09);
\shade[shading=ball, ball color=black] (3,-1.73) circle (.09);
\shade[shading=ball, ball color=black] (5,-1.73) circle (.09);
\shade[shading=ball, ball color=black] (4,3.46) circle (.09);
\shade[shading=ball, ball color=black] (0,0) circle (.09);
\shade[shading=ball, ball color=black] (2,0) circle (.09);
\shade[shading=ball, ball color=black] (4,0) circle (.09);
\draw [fill = yellow] (0,0) -- (2,0) -- (1,1.73) -- (0,0);
\draw [fill = yellow] (3,1.73) -- (2,0) -- (1,1.73) -- (3,1.73);
\draw [fill = yellow] (3,1.73) -- (4,3.46) -- (5,1.73) -- (3,1.73);
\draw [fill = yellow] (4,0) -- (3,1.73) -- (5,1.73) -- (4,0);
\draw [fill = yellow] (2,0) -- (4,0) -- (3,-1.73) -- (2,0);
\draw [fill = yellow] (4,0) -- (3,-1.73) -- (5,-1.73) -- (4,0);
\end{tikzpicture}
\end{center}

Now it is easy to check that $ii''(i+1)(i+1)'  \in L(S)$ for every $i \in \Z_3$. It follows that
$123$ is a transversal of the successive differences for the chain
$$\emptyset \subset 1 \subset 11''22' \subset V$$
in $L(S)$. 

Finally, each facet of the form $ii''(i+1)p$ or $ii''(i+1)'p$ is a transversal of the successive differences for the chain
$$\emptyset \subset i \subset ii'' \subset ii''(i+1)(i+1)' \subset V$$
in $L(S)$. Since we have now checked all facets, it follows that $S$ is a BRSC.

Let $\cl(X)$ denote the closure of $X \subseteq V$ in $L(T_3(S))$. For each $i \in \Z_3$, we have $i(i+1)(i+1)', i''(i+1)(i+1)' \notin \H$, so we successively get $(i+1)' \in \cl(i(i+1))$ and $i'' \in \cl(i(i+1))$. Thus $\cl(i(i+1))$ contains $ii''(i+1) \in \fct(T_3(S))$, yielding $\cl(i(i+1)) = V$. But then $i \in \cl(123 \setminus \{ i \})$ for every $i \in 123$. Since $123 \in T_3(\H)$, then there is no chain of the form (\ref{derby3}) and so $T_3(S)$ is not boolean representable.

We remark that $S$ is not pure since it is straightforward to check that $1'2'2''$ is a facet. But $T_3(S)$ is pure because there are no facets of dimension 1: given distinct $p,q \in V$, there exists some $r \in V \setminus pq$ such that $pqr$ is not a yellow triangle.
\qed

\noindent
{\em Analysis of Example \ref{bfour}}.
Let $M$ be such a boolean matrix. Since the columns are all distinct and nonzero, every pair of distinct columns is independent. Now let $X$ be a set of independent columns with $|X| = 2$ or $3$. Let $I \subset 1234$ be such that the square matrix $M[I,X]$ is nonsingular. Let $j \in 1234 \setminus I$ and let $c$ be the column of $M$ having a 1 at row $j$ and 0 elsewhere. Then the permanent of $M[I \cup \{ j \},X \cup \{ c \}]$ equals the permanent of $M[I,X]$, hence $M[I \cup \{ j \},X \cup \{ c \}]$ is nonsingular. and so  $X \cup \{ c \}$ is independent. Thus $\B^{(4)}$ is pure.

Since $\B^{(4)}$ is by definition a BRSC, then $\B^{(4)}_3$ is a TBRSC. Let $\oo{X}$ denote the closure of $X$ in $\flatx \B^{(4)}_3$. Consider the columns of $M$ defined by
$$a = \begin{bmatrix}
1\\ 0\\ 0\\ 0
\end{bmatrix}, \quad
b = \begin{bmatrix}
1\\ 1\\ 1\\ 0
\end{bmatrix}, \quad
c = \begin{bmatrix}
1\\ 1\\ 0\\ 1
\end{bmatrix}.$$
The permanent of the matrix
$$M[134,abc] = \begin{bmatrix}
1 & 1 & 1\\ 0 & 1 & 0\\ 0 & 0 & 1
\end{bmatrix}$$
is 1, hence $abc$ is independent. Define
$$d = \begin{bmatrix}
0\\ 1\\ 1\\ 0
\end{bmatrix}, \quad
e = \begin{bmatrix}
1\\ 0\\ 1\\ 0
\end{bmatrix}, \quad
f = \begin{bmatrix}
0\\ 0\\ 1\\ 1
\end{bmatrix} \quad
g = \begin{bmatrix}
1\\ 0\\ 1\\ 1
\end{bmatrix}.$$
We have
$$M[1234,abd] = \begin{bmatrix}
1 & 1 & 0\\ 0 & 1 & 1\\ 0 & 1 & 1\\ 0 & 0 & 0
\end{bmatrix},\quad 
M[1234,bde] = \begin{bmatrix}
1 & 0 & 1\\ 1 & 1 & 0\\ 1 & 1 & 1\\ 0 & 0 & 0
\end{bmatrix},\quad 
M[1234,abe] = \begin{bmatrix}
1 & 1 & 1\\ 0 & 1 & 0\\ 0 & 1 & 1\\ 0 & 0 & 0
\end{bmatrix}.$$
Since no row of $M[1234,abd]$ has precisely two zeroes, $abd$ is dependent. The same occurs with $bde$. It is immediate that $M[123,abe]$ has permanent 1, hence $abe$ is independent. Thus we successively deduce $d \in \oo{ab}$, $e \in \oo{ab}$ and so $\oo{ab}$ contains the facet $abe$. Therefore $\oo{ab} = V$, where $V$ denotes the full set of vertices. Out of symmetry, so is $\oo{ac}$. 

Now
$$M[1234,bcf] = \begin{bmatrix}
1 & 1 & 0\\ 1 & 1 & 0\\ 1 & 0 & 1\\ 0 & 1 & 1
\end{bmatrix},\quad 
M[1234,bcg] = \begin{bmatrix}
1 & 1 & 1\\ 1 & 1 & 0\\ 1 & 0 & 1\\ 0 & 1 & 1
\end{bmatrix},\quad 
M[1234,bfg] = \begin{bmatrix}
1 & 0 & 1\\ 1 & 0 & 0\\ 1 & 1 & 1\\ 0 & 1 & 1
\end{bmatrix}.$$
Since no row of $M[1234,bcf]$ has precisely two zeroes, $bcf$ is dependent. The same occurs with $bcg$. It is immediate that $M[123,bfg]$ has permanent 1, hence $bfg$ is independent. Thus we successively deduce $f \in \oo{bc}$, $g \in \oo{bc}$ and so $\oo{bc}$ contains the facet $bfg$. Therefore $\oo{bc} = V$. Together with $\oo{ab} = \oo{ac} = V$, and $abc$ being independent, this proves that $\B^{(4)}_3$ is not a BRSC.
\qed

\noindent
{\em Proof of Theorem \ref{spch}}.
Suppose first that $k \leq 2$. By Proposition \ref{tbone}, $T_k(S)$ is a BRSC, therefore {\rm pure}$(T_k(S))$ is a BRSC by Corollary \ref{kan}(ii).

Thus we may assume that $k = 3$. Write {\rm pure}$(T_3(S)) = (V',\H')$ and consider the restriction $S|_{V'}$. Then  {\rm pure}$(T_3(S)) = {\rm pure}(T_3(S|_{V'})_3)$. Since 
BRSCs are closed under restriction, $S|_{V'}$ is also a BRSC. Therefore we may assume that $V' = V$.

Let 
$$\F = F \in \L(S) \; {\big{\lvert}} \; |F \cap X| \neq 1 \mbox{ for every }X \in \fct(S) \cap P_2(V) \}.$$
We claim that $\F$ is a Moore family.

Clearly, $V \in \F$. Let $F,F' \in \F$. We have $F \cap F' \in L(S)$. Let $X \in \fct(S) \cap P_2(V)$.  Suppose that $|(F \cap F') \cap X| = 1$. Then $|F \cap X| = 1$ or $|F' \cap X| = 1$, contradicting $F,F' \in \F$. Thus $F \cap F' \in \F$ and so $\F$ is a Moore family.

Therefore $S' = (V,{\rm Tr}(\F))$ is a BRSC. We claim that $\pure(T_3(S)) = T_3(S')$.

For every $Y \subseteq V$, let $\oo{Y}$ denote its closure in $L(S)$. Let $X \in \H \cap P_3(V)$. Then there exists an enumeration $a,b,c$ of the elements of $X$ such that 
\beq
\label{lenha1}
\emptyset \subset \oo{a} \subset \oo{ab} \subset V
\eeq
and $c \notin \oo{ab}$. Clearly, $\emptyset,V \in \F$. 

Suppose that $Y \in \fct(S) \cap P_2(V)$ satisfies $|Y \cap \oo{a}| = 1$. We may write $Y = yz$ with $y \in \oo{a}$. If $z \in \oo{ab}$ (respectively $z \notin \oo{ab}$), then $yzc$ (respectively $ybz$) is a transversal of the successive differences for (\ref{lenha1}), contradicting $yz \in \fct(S)$. Thus $\oo{a} \in \F$.

Suppose now that $Y \in \fct(S) \cap P_2(V)$ satisfies $|Y \cap \oo{ab}| = 1$. We may write $Y = yz$ with $y \in \oo{ab}$. If $y \in \oo{a}$ (respectively $y \notin \oo{a}$), then $ybz$ (respectively $ayz$) is a transversal of the successive differences for (\ref{lenha1}), contradicting $yz \in \fct(S)$. Thus $\oo{ab} \in \F$.

Therefore (\ref{lenha1}) is a chain in $\F$ and so $X \in {\rm Tr}(\F)$. It follows that $\H' \subseteq {\rm Tr}(F) \cap P_{\leq 3}(V)$.

Conversely, let $X \in {\rm Tr}(F) \cap P_{\leq 3}(V)$. Since $\F \subseteq L(S)$, we have ${\rm Tr}(\F) \subseteq \H$ and so $X \in T_3(\H)$. We certainly have $X \in \H'$ if if $|X| = 0$ or 3, and the case $|X| = 1$ follows from $V' = V$. Hence we may assume that $|X| = 2$. There exists an enumeration $a,b$ of the elements of $X$ and $F \in \F$ such that $a \in F$ and $b \notin F$. But then, by definition of $\F$, we get $X \notin \fct(S)$. Hence there exists some $c \in V \setminus X$ such that $X \cup \{ c \} \in \H$. Thus $X \cup \{ c \} \in \H'$ and so $X \in \H'$. Therefore $T_3({\rm Tr}(F)) \subseteq \H'$ and so $\pure(T_3(S)) = T_3(S')$. It follows that $\pure(T_3(S))$ is a TBRSC.
\qed

\noindent
{\em Proof of Lemma \ref{connected}}.
Obviously, $S$ is disconnected if $\H = P_1(V)$ and $|V| > 1$, and
connected if $|V| = 1$. Hence
we may assume that $pq \in \H$ for some distinct $p,q \in V$. 

Let $M$ be an $R \times V$ boolean matrix representing $S$. It
follows from $pq \in \H$ that $M[R,p] \neq M[R,q]$. Thus, for every $v
\in V$, we have either $M[R,v] \neq M[R,p]$ or $M[R,v] \neq M[R,q]$,
implying that $vp$ or $vq$ is an edge in $\H$.
Therefore $S$ is connected.
\qed

\noindent
{\em Proof of Theorem \ref{wed}}.
(i) It follows from Lemma \ref{connected} that $S$ is connected.
By Hurewicz Theorem (see \cite{Hat}), the 1st homology group of
$S$ is the abelianization of $\pi_1(||S||)$, and therefore, in view
of Theorem \ref{fung}, a free abelian group of known rank. The second
homology group of any 2-dimensional simplicial complex is
$\ker\partial_2 \leq C_2(S)$, that is, a subgroup of a free abelian
group. Therefore $H_2(S)$ is itself free abelian.

(ii) By \cite[Proposition 3.3]{Wal}, any finite 2-dimensional
simplicial complex with free fundamental group has the homotopy type
of a wedge of 1-spheres and 2-spheres.
\qed

\vspace{1cm}

{\sc Stuart Margolis, Department of Mathematics, Bar Ilan University,
  52900 Ramat Gan, Israel}

{\em E-mail address:} margolis@math.biu.ac.il

\bigskip

{\sc John Rhodes, Department of Mathematics, University of California,
  Berkeley, California 94720, U.S.A.}

{\em E-mail addresses}: rhodes@math.berkeley.edu, BlvdBastille@gmail.com

\bigskip

{\sc Pedro V. Silva, Centro de
Matem\'{a}tica, Faculdade de Ci\^{e}ncias, Universidade do
Porto, R. Campo Alegre 687, 4169-007 Porto, Portugal}

{\em E-mail address}: pvsilva@fc.up.pt

\end{document}